\documentclass{article}

\usepackage{amsmath, amsthm, amsfonts, amssymb}
\usepackage[all]{xy}

\theoremstyle{plain}
\newtheorem{thm}{Theorem}[section]
\newtheorem{lem}[thm]{Lemma}
\newtheorem{prop}[thm]{Proposition}
\newtheorem{cor}[thm]{Corollary}
\newtheorem*{question}{Question}

\theoremstyle{definition}
\newtheorem*{defn}{Definition}

\newtheoremstyle{citing}
	{3pt}{3pt}{\itshape}{}{\bfseries}{.}{.5em}{\thmnote{#3}}

\theoremstyle{citing}
\newtheorem*{thmx}{}

\newcommand{\sbinom}[2]{\genfrac{[}{]}{0pt}{}{#1}{#2}}

\newcommand{\rank}{\mathrm{rank}}
\newcommand{\Aut}{\mathrm{Aut}}
\newcommand{\Hom}{\mathrm{Hom}}
\newcommand{\GL}{\mathrm{GL}}
\newcommand{\Der}{\mathrm{Der}}
\newcommand{\emb}{\mathrm{emb}}
\newcommand{\qemb}{\mathrm{qemb}}
\newcommand{\pow}{\mathrm{pow}}
\newcommand{\Fcom}{\mathrm{Fcom}}
\newcommand{\com}{\mathrm{com}}
\newcommand{\proj}{\mathrm{proj}}
\newcommand{\F}{\mathbb{F}}

\newcommand{\Z}{\mathbb{Z}}
\newcommand{\ep}{\varepsilon}
\allowdisplaybreaks[1]

\title{The Automorphism Group of a Finite $p$-Group is Almost Always a $p$-Group}
\author{Geir T. Helleloid \\
				Department of Mathematics, Bldg. 380 \\
				Stanford University \\
				Stanford, CA 94305-2125 \\
				\texttt{geir@math.stanford.edu}
\and
				Ursula Martin \\
				Department of Computer Science \\
				Queen Mary University of London \\
				Mile End Road \\
				London E1 4NS, UK \\
				\texttt{Ursula.Martin@dcs.qmul.ac.uk}
}

\begin{document}
\maketitle

\begin{abstract}
Many common finite $p$-groups admit automorphisms of order coprime to $p$, and when $p$ is odd, it is reasonably difficult to find finite $p$-groups whose automorphism group is a $p$-group.  Yet the goal of this paper is to prove that the automorphism group of a finite $p$-group is almost always a $p$-group.  The asymptotics in our theorem involve fixing any two of the following parameters and letting the third go to infinity: the lower $p$-length, the number of generators, and $p$.  The proof of this theorem depends on a variety of topics: counting subgroups of a $p$-group; analyzing the lower $p$-series of a free group via its connection with the free Lie algebra; counting submodules of a module via Hall polynomials; and using numerical estimates on Gaussian coefficients.
\end{abstract}

\section{Introduction}
\label{intro}

The goal of this paper is to prove that, in a certain asymptotic sense, the automorphism group of a finite $p$-group is almost always a $p$-group.  A weaker version of this result was announced by the second author in \cite{mar}, but this paper contains the first published proof.

The result may not seem entirely plausible at first, as many common finite $p$-groups have an automorphism group that is not a $p$-group.  Examples include: abelian $p$-groups, unless $p = 2$ and the type of the group has repeated parts (see Macdonald~\cite[Chapter II, Theorem 1.6]{mac}); the Sylow-$p$ subgroup of $\GL(n, \F_p)$ for $p$ odd (see Gibbs~\cite{gib}); and the extraspecial $p$-groups (see Winter~\cite{win}).  Furthermore, Bryant and Kov\'{a}cs~\cite{bk} show that any finite group occurs as the quotient $A(H)$ of the automorphism group of some finite $p$-group $H$, where $A(H)$ is as defined below.  Our result seems to say that most $p$-groups are complicated and unnatural-looking and that familiar examples are far from typical.

It is reasonably easy to find finite $2$-groups whose automorphism group is a $2$-group: $\Z_{2^n}$, the dihedral 2-group $D_{2^n}$ ($n \ge 3$), and the generalized quaternion group $Q_{2^n}$ ($n \ge 4$) are common examples, while Newman and O'Brien~\cite{no} offer three more infinite families.  It is more difficult to find finite $p$-groups whose automorphism groups are $p$-groups when $p$ is odd.  In~\cite{hor}, Horo{\v{s}}evski{\u\i} constructs such a $p$-group with nilpotence class $n$ for each $n \ge 2$ and such a $p$-group on $d$ generators for each $d \ge 3$.  Furthermore, Horo{\v{s}}evski{\u\i} shows in \cite{hor} and \cite{hor2} that for any prime $p$, if $H_1, H_2, \dots, H_n$ are finite $p$-groups whose automorphism groups are $p$-groups, then the automorphism group of the iterated wreath product $H_1 \wr H_2 \wr \cdots \wr H_n$ is also a $p$-group.  Otherwise, most known examples arise from complicated and unnatural-looking constructions (see Webb~\cite{web}).  A survey on the automorphism groups of finite $p$-groups, including a comprehensive list of examples in the literature of finite $p$-groups whose automorphism groups are $p$-groups, can be found in~\cite{hel}.

In a computational vein, Eick, Leedham-Green, and O'Brien~\cite{elo} describe an algorithm for constructing the automorphism group of a finite $p$-group.  This algorithm has been implemented by Eick and O'Brien in the GAP package AutPGroup~\cite{GAP4}.  Compiled with the gracious help of Eamonn O'Brien (personal communication) and the GAP packages AutPGroup and SmallGroups~\cite{GAP4}, Table~\ref{data} summarizes data on the proportion of small $p$-groups whose automorphism group is a $p$-group.  (More information about the SmallGroups package can be found in Besche, Eich and O'Brien~\cite{beo}.)

\begin{table} \label{data}
\begin{center}
\begin{tabular}{|c|c|c|c|}
\hline
Order & $p = 2$ & $p = 3$ & $p = 5$ \\
\hline
$p^3$ & 3 of 5 & 0 of 5 & 0 of 5 \\
$p^4$ & 9 of 14 & 0 of 15 & 0 of 15 \\
$p^5$ & 36 of 51 & 0 of 67 & 1 of 77 \\
$p^6$ & 211 of 267 & 30 of 504 & 65 of 685 \\
$p^7$ & 2067 of 2328 & 2119 of 9310 & 11895 of 34297 \\
\hline
\end{tabular}
\parbox{3in}{\caption{The proportion of $p$-groups of a given order whose automorphism group is a $p$-group.}}
\end{center}
\end{table}

Of course, the meaning of the statement ``The automorphism group of a finite $p$-group is almost always a $p$-group'' depends on the asymptotic interpretation of ``almost always.''  Probably the most natural interpretation is to consider all $p$-groups of order at most $p^n$ and let $n$ go to infinity.  However, this is not the sense of our result, and indeed, the question remains open for this interpretation (see Mann~\cite[Question 9]{man}).  The precise statement of our main theorem depends on the \emph{lower $p$-series} of a group.  The lower $p$-series will be defined in Section~\ref{lowerpsec}; for the moment, it suffices to say that the lower $p$-series is a central series with elementary abelian factors and that the \emph{lower $p$-length} of a group is the number of non-identity terms in the associated lower $p$-series.  The main theorem of this paper may be concisely stated as follows.

\begin{thm} \label{mainthm}
Fix a prime $p$ and positive integers $d$ and $n$.  Let $r_{d, n}$ be the proportion of $p$-groups minimally generated by $d$ elements and with lower $p$-length at most $n$ whose automorphism group is a $p$-group.  If $n \ge 2$, then
\[
\lim_{d \to \infty}{r_{d, n}} = 1.
\]
If $d \ge 5$, then
\[
\lim_{n \to \infty}{r_{d,n}} = 1.
\]
If $n = 2$ and $d \ge 10$, or $n \ge 3$ and $d \ge 6$, or $n \ge 10$ and $d \ge 5$, then
\[
\lim_{p \to \infty}{r_{d,n}} = 1.
\]
\end{thm}

The proof of Theorem~\ref{mainthm} breaks down into three parts, which are presented in Sections~\ref{lowerpsec}, \ref{est1sec}, and \ref{est2sec}, and are assembled to prove Theorem~\ref{mainthm} in Section~\ref{summarysec}.  In the remainder of this section, we will outline the structure of the proof.

The first step is to connect the enumeration of finite $p$-groups to an analysis of certain subgroups and quotients of free groups.  Let $F$ be the free group on $d$ generators and let $F_n$ be the $n$-th term in the lower $p$-series of $F$.  It turns out that the action of $\Aut(F/F_{n+1})$ on $F_n/F_{n+1}$ induces an action of $\GL(d, \F_p)$ on $F_n/F_{n+1}$, and the $\Aut(F/F_{n+1})$-orbits on the normal subgroups of $F_n/F_{n+1}$ are also the $\GL(d, \F_p)$-orbits.

For any finite $p$-group $H$, write $A(H)$ for the group of automorphisms of $H/\Phi(H)$ induced by $\Aut(H)$, where $\Phi(H)$ is the Frattini subgroup of $H$.  We shall see that if $A(H)$ is a $p$-group then so is $\Aut(H)$; in fact, our main goal is to prove, in some sense, that $A(H)$ is usually trivial.  In Section~\ref{lowerpsec}, after defining and investigating the lower $p$-series, we prove the following theorem.

\begin{thm} \label{bijthm}
Fix a prime $p$ and integers $d, n \ge 2$.  Let $F$ be the free group on $d$ generators and define the following sets:
\begin{eqnarray*}
\mathcal{A}_{d,n} &=& \{\textrm{normal subgroups of $F/F_{n+1}$ lying in $F_2/F_{n+1}$}\} \\
\mathcal{B}_{d,n} &=& \{\textrm{normal subgroups of $F/F_{n+1}$ lying in $F_2/F_{n+1}$} \\
&& \qquad \textrm{and not containing $F_n/F_{n+1}$}\} \\
\mathcal{C}_{d,n} &=& \{\textrm{normal subgroups of $F/F_{n+1}$ lying in $F_n/F_{n+1}$}\} \\
\mathcal{D}_{d,n} &=& \{\textrm{normal subgroups of $F/F_{n+1}$ contained in the} \\
&& \qquad \textrm{regular $\GL(d,\F_p)$-orbits in $\mathfrak{C}_{d,n}$}\} \\
\\
\mathfrak{A}_{d,n} &=& \{\textrm{$\Aut(F/F_{n+1})$-orbits in $\mathcal{A}_{d,n}$}\} \\
\mathfrak{B}_{d,n} &=& \{\textrm{$\Aut(F/F_{n+1})$-orbits in $\mathcal{B}_{d,n}$}\} \\
\mathfrak{C}_{d,n} &=& \{\textrm{$\Aut(F/F_{n+1})$-orbits in $\mathcal{C}_{d,n}$}\} = \{\textrm{$\GL(d,\F_p)$-orbits in $\mathcal{C}_{d,n}$}\} \\
\mathfrak{D}_{d,n} &=& \{\textrm{regular $\GL(d,\F_p)$-orbits in $\mathcal{C}_{d,n}$}\}.
\end{eqnarray*}
Then there is a well-defined map $\pi_{d,n}: \mathfrak{A}_{d,n} \to \{\textrm{finite $p$-groups}\}$ given by $L/F_{n+1} \mapsto F/L$, where $L/F_{n+1} \in \mathcal{A}_{d,n}$.  Furthermore $\pi_{d,n}$ induces bijections
\begin{eqnarray*}
\mathfrak{A}_{d,n} &\leftrightarrow& \{\textrm{$p$-groups of lower $p$-length at most $n$} \\
&& \qquad \textrm{and minimally generated by $d$ elements} \} \\
\mathfrak{B}_{d,n} &\leftrightarrow& \{\textrm{$p$-groups of lower $p$-length $n$} \\
&& \qquad \textrm{and minimally generated by $d$ elements} \} \\
\mathfrak{D}_{d,n} &\leftrightarrow& \{\textrm{subgroups $H$ in $\pi_{d,n}(\mathfrak{C}_{d,n})$ with $A(H) = 1$}\}.
\end{eqnarray*}
\end{thm}
Recall that a regular orbit is one in which every point has trivial stabilizer.  Note that as a result of Theorem~\ref{bijthm}, it will be enough to show that $|\mathfrak{A}_{d,n}|/|\mathfrak{D}_{d,n}|$ goes to 1 under the relevant limits.

Section~\ref{freegpsec} follows with an examination of the structure of $F_n/F_{n+1}$ that will be needed in Section~\ref{est1sec}.  Section~\ref{est0sec} contains combinatorial estimates, including bounds on Gaussian coefficients, that will be needed in Sections~\ref{est1sec} and~\ref{est2sec}.  Finally, the second and third steps of the proof of Theorem~\ref{mainthm} are summarized in Theorems~\ref{est1thm} and~\ref{est2thm} and are proved in Sections~\ref{est1sec} and~\ref{est2sec}.  The terms $C(p)$ and $D(p)$ that appear in Theorems~\ref{est1thm} and~\ref{est2thm} are functions of $p$ which tend to 1 as $p \to \infty$.

\begin{thm} \label{est1thm}
Fix a prime $p$ and integers $d$ and $n$ so that either $n \ge 3$ and $d \ge 6$ or $n \ge 10$ and $d \ge 5$.  Let $F$ be the free group on $d$ generators and let $d_n$ be the rank of $F_n/F_{n+1}$.  Then
\[
1 \le \frac{|\mathfrak{A}_{d,n}|}{|\mathfrak{C}_{d,n}|} \le 1 + C(p)^{n-1} D(p)^{n-2} p^{d_{n-1} - d_n/4 + d^2}.
\]
\end{thm}

The proof of Theorem~\ref{est1thm} uses a theorem estimating the number of normal subgroups of an arbitrary finite $p$-group, applying it to quotients of free groups.

\renewcommand{\arraystretch}{1.2}
\begin{thm} \label{est2thm}
Fix a prime $p$ and integers $d$ and $n$ so that either $n = 2$ and $d \ge 10$ or $n \ge 3$ and $d \ge 3$.  Let $F$ be the free group on $d$ generators and let $d_n$ be the rank of $F_n/F_{n+1}$.  Let
\[
K = \left\{
	\begin{array}{r@{\quad:\quad}l}
		C(p)^5 D(p)^4 p^{17/4} & \textrm{$n = 2$ and $d \ge 10$} \\
		C(p)^2 D(p) p^{3/4} & n \ge 3.
	\end{array}
	\right.
\]
Let
\[
x = \left\{
	\begin{array}{r@{\quad:\quad}l}
		-d & n = 2 \\
		d^2-d_n/2 & n \ge 3.
	\end{array}
	\right.
\]
Then
\begin{enumerate}
\renewcommand{\labelenumi}{\emph{(\alph{enumi})}}
\renewcommand{\theenumi}{\ref{est2thm}(\alph{enumi})}
\item
\[
1 \le \frac{|\mathfrak{C}_{d, n}| \cdot |\GL(d,\F_p)|}{|\mathcal{C}_{d, n}|} \le 1 + K p^x.
\]
\item
\[
1 \le \frac{|\mathfrak{C}_{d, n}|}{|\mathfrak{D}_{d, n}|} \le \frac{1+K p^x}{1-K p^x}.
\]
\end{enumerate}
\end{thm}
\renewcommand{\arraystretch}{1}

In stating Theorems~\ref{est1thm} and~\ref{est2thm}, we have judged it more satisfactory to give explicit numerical bounds, even though the proof of Theorem~\ref{mainthm} requires only asymptotic bounds.  However, since we have no expectation that our proof method gives bounds that are sharp, we have opted for clean explicit bounds rather than the best possible.

As we will show in Section~\ref{summarysec}, Theorem~\ref{mainthm} follows easily from Theorems~\ref{bijthm}, \ref{est1thm}, and~\ref{est2thm}.  We close Section~\ref{summarysec} with some observations and open questions.

\section{The Lower $p$-Series}
\label{lowerpsec}

In this section, we define and discuss the \emph{lower $p$-series} of a group (also called the \emph{lower central $p$-series} or the \emph{lower exponent-$p$ central series}).  Then, in Theorems~\ref{isoenum1} and \ref{isoenum2}, we describe how isomorphism classes of finite $p$-groups in a variety may be enumerated, obtaining Theorem~\ref{bijthm} as a corollary.  

\subsection{Preliminaries}

The lower $p$-series was introduced by Skopin~\cite{sko} and Lazard~\cite{laz}, and it is described in detail by Huppert and Blackburn~\cite[Chapter VIII]{hb} (under the name $\lambda$-series) and by Bryant and Kov\'{a}cs~\cite{bk}.  The lower $p$-series is particularly suited to computer analysis of finite $p$-groups and forms the basis of the $p$-group generation algorithm of M. F. Newman~\cite{new} (this algorithm is described in greater detail in, for example, O'Brien~\cite{obr}).  This algorithm was modified in~\cite{obr2} and~\cite{elo} to construct automorphism groups of finite $p$-groups.  It should also be mentioned that some information about the lower $p$-series has appeared in~\cite{obr} and~\cite{elo}, while the link between the lower $p$-series and automorphisms described in Subsection~\ref{autgroup} is an extension of results that Higman~\cite{hig} and Sims~\cite{sim} used to count finite $p$-groups.

\begin{defn}
Fix a prime $p$.  For any group $H$, the \emph{lower $p$-series} $H = H_1 \ge H_2 \ge \cdots$ of $H$ is defined by $H_{i+1} = H_i^p [H_i, H]$ for $i \ge 1$.  $H$ is said to have \emph{lower $p$-length} $n$ if $H_n$ is the last non-identity element of the lower $p$-series.
\end{defn}

Note that if $H$ is a finite $p$-group, then $H_2 = \Phi(H)$, the Frattini subgroup of $H$.  Before we list some basic facts about the lower $p$-series, recall that a subgroup is \emph{fully invariant} if every endomorphism of the group restricts to an endomorphism of the subgroup.  Also, we will write $H = \gamma_1(H) \ge \gamma_2(H) \ge \cdots$ to denote the \emph{lower central series} of $H$, where $\gamma_{i+1}(H) = [\gamma_i(H), H]$.  The following proposition states five fundamental properties of the lower $p$-series; the first four facts are proved in Huppert and Blackburn~\cite[Chapter VIII, Theorem 1.5 and Corollary 1.6]{hb} and the fifth fact is obvious by induction.

\begin{prop} \label{fratprops}
For all positive integers $i$ and $j$,
\begin{enumerate}
\item $[H_i, H_j] \le H_{i+j}$.
\item $H_i^{p^j} \le H_{i+j}$.
\item $H_i = \gamma_1(H)^{p^{i-1}} \gamma_2(H)^{p^{i-2}} \cdots \gamma_i(H)$.
\item $H_{i+1}$ is the smallest normal subgroup of $H$ lying in $H_i$ such that $H_i/H_{i+1}$ is an elementary abelian $p$-group and is central in $H/H_{i+1}$.
\item $H_i$ is fully invariant in $H$.
\end{enumerate}
\end{prop}

As we will see, the fact that $H_i/H_{i+1}$ is elementary abelian, and therefore an $\F_p$-vector space, is a key reason we are able to prove the main theorem.  It is easy to see the following proposition.

\begin{prop}
Let $H$ be a finite group.  Then $H$ is a $p$-group if and only if $H$ has finite lower $p$-length.
\end{prop}

The lower $p$-length of a finite $p$-group is related to the lower $p$-series of a free group in the following way.  Let $F$ be the free group on $d$ generators; then any finite $p$-group $H$ that is minimally $d$-generated is isomorphic to $F/U$ for some normal subgroup $U$ of $F$.  By induction, $H_i = F_iU/U$:
\begin{eqnarray*}
H_{i+1} &=& (F_iU/U)^p [F_iU/U, F/U] \\
&=& F_i^p [F_i, F] U/U \\
&=& F_{i+1}U/U.
\end{eqnarray*}
So the lower $p$-length of $H$ is $n$, where $F_{n+1}$ is the first term in the lower $p$-series of $F$ that is contained in $U$.

\subsection{The Lower $p$-Series and Automorphisms}
\label{autgroup}

In this subsection we collect some necessary facts linking the lower $p$-series and automorphisms.  First, suppose that $H$ is a finite $p$-group that is minimally $d$-generated.  Of course, every automorphism of $H$ induces an automorphism of $H_i/H_{i+1}$ for each $i$.  In particular, any automorphism of $H$ induces an element of $\Aut(H/H_2) \cong \GL(d, \F_p)$ (by the Burnside Basis Theorem, the rank of $H/H_2$ is $d$).  Thus we obtain a map from $\Aut(H)$ to $\GL(d, \F_p)$, and an exact sequence
\[
1 \to K(H) \to \Aut(H) \to A(H) \to 1,
\]
where $A(H)$ is a subgroup of $\GL(d, \F_p)$.  The group $K(H)$ acts trivially on $H/H_2$, and hence on each factor $H_i/H_{i+1}$ (see Huppert and Blackburn~\cite[Chapter VIII, Theorem 1.7]{hb}).  As $\Aut(H)$ acts on each $H_i/H_{i+1}$ and the kernel of the action contains $K(H)$, we obtain an action of $A(H)$ on each $H_i/H_{i+1}$.  The following key proposition is due to P. Hall~\cite[Section 1.3]{hal}.

\begin{prop}
If $H$ is a finite $p$-group, then so is $K(H)$.
\end{prop}

Let $F$ be the free group on $d$ generators $y_1, y_2, \dots, y_d$.  We need two observations about the subgroup $F_2$, first recalling an obvious result on the Frattini quotient.

\begin{prop}
\label{frataut}
If $H$ is a finite $p$-group and $\theta$ is an endomorphism of $H$ that induces an automorphism on the Frattini quotient $H/H_2$, then $\theta$ is an automorphism of $H$.
\end{prop}

\begin{prop}
\label{f2maxl}
$F_2$ is a maximal fully invariant subgroup of $F$.
\end{prop}

\begin{proof}
Suppose $U > F_2$ is a fully invariant subgroup of $F$.  The elements $y_1^{a_1} \cdots y_d^{a_d}$, with $0 \le a_i < p$, form a complete set of coset representatives for the cosets of $F_2$ in $F$, so $U$ contains an element $y = y_1^{a_1} \cdots y_d^{a_d}$ with some $a_i$ nonzero.  Fix $1 \le k \le d$ and let $b_i$ be a multiplicative inverse of $a_i$ modulo $p$.  Then the endomorphism of $F$ that sends $y_j$ to 1 for $j \neq i$ and sends $y_i$ to $y_k^{b_i}$ also sends $y$ to $y_k$, showing that $y_k \in U$.  This holds for $1 \le k \le d$, and so $U = F$.
\end{proof}

\begin{prop}
\label{f2lift}
Let $U$ be a fully invariant subgroup of $F$ contained in $F_2$ with $H = F/U$ a finite $p$-group.  Then any automorphism $\theta$ of $F/F_2$ lifts to an automorphism of $H$.
\end{prop}

\begin{proof}
Since $F$ is free, there is an endomorphism $\theta'$ of $F$ such that $\theta'(y_i) \in \theta(y_i F_2)$ for $1 \le i \le d$.  Therefore $\theta'(y) \in \theta(y F_2)$ for all $y \in F$.  Then $\theta'$ induces $\theta$ on $F/F_2$, and since $U$ is fully invariant, maps $U$ to itself.  So $\theta'$ induces an endomorphism $\theta''$ of $H$.  But $\theta''$ induces $\theta$, an automorphism of $F/F_2 \cong (F/U)/(F_2/U) \cong H/H_2$, the Frattini quotient of $H$.  By Proposition~\ref{frataut}, $\theta''$ is an automorphism of $H$.  Thus $\theta$ lifts to an automorphism $\theta''$ of $H$.
\end{proof}

Finally, we note that by Huppert and Blackburn~\cite[Chapter VIII, Theorem 11.15]{hb}, the rank of $F/[F,F]$, and hence of $F/F_2$, is $d$, and the rank of $F_n/F_{n+1}$ is finite for each $n$ (in Section~\ref{freegpsec}, we will compute the rank of $F_n/F_{n+1}$ in general).

\subsection{Enumerating Groups in a Variety}

A \emph{variety of groups} $V$ consists of all groups $G$ satisfying a set of relations $w = 1$, where $w$ ranges over a fixed set $W$ of group words (see Neumann~\cite{neu}).  Let $F$ be the free group on $d$ generators.  The variety $V$ contains a \emph{relatively free group} on $d$ generators, namely $F/U$, where $U$ is the \emph{verbal subgroup} of $F$ generated by all the values of $w \in W$.  For example, all abelian groups form the variety in which the relation $ab = ba$ holds for all group elements $a$ and $b$.  Then the free abelian group on $d$ generators is the relatively free group on $d$ generators in the variety of abelian groups.  We will only be interested in the variety of $p$-groups of lower $p$-length at most $n$, but the theorems in this subsection hold in more general situations.

Let $U$ be a fully invariant subgroup of $F$.  Then $G = F/U$ is a relatively free group in some variety $V$ on at most $d$ generators.  The relations defining $V$ come from setting each word in $U$ equal to the identity element.  Suppose that $G$ is a finite non-trivial $p$-group.  In this setting, we can describe $A(G)$ and $K(G)$ more precisely.

Note that $F_2 U$ is a fully invariant subgroup of $F$, and by Proposition~\ref{f2maxl}, either $F = F_2 U$ or $F_2 = F_2 U$.  In the first case, $F = U$, contradicting the non-triviality of $G$.  Thus $F_2 = F_2 U$ and $U \le F_2$.  Since $F/F_2$ has rank $d$, both $F$ and $F/F_2$ are minimally generated by $d$ elements.  It follows that $G = F/U$ is also minimally generated by $d$ elements.

\begin{thm} \label{isoenum1}
Suppose that $G$ is the relatively free group on $d$ generators in a variety of groups $V$ and that $|G| = p^g$.  Then 
\[
1 \to K(G) \to \Aut(G) \to \GL(d,\F_p) \to 1
\]
is exact and $|K(G)| = p^{d(g-d)}$.  Furthermore, the map $L \mapsto G/L$ defines a bijection between $\Aut(G)$-orbits of normal subgroups $L$ of $G$ lying in $G_2$ and $d$-generator groups in $V$.  If $H = G/L$, then 
\[
1 \to B(L) \to N_{\Aut(G)}(L) \to \Aut(H) \to 1
\]
is exact, where $B(L)$ is the subgroup of $N_{\Aut(G)}(L)$ that acts trivially on $H$.  If $|L| = p^m$, then $|B(L)| = p^{dm}$.
\end{thm}

\begin{proof}
By Proposition~\ref{f2lift}, any automorphism $\theta$ of $F/F_2 \cong G/G_2$ lifts to an automorphism of $G$.  Thus $A(G)$ is the full automorphism group of $F/F_2$, which is $\GL(d,\F_p)$.  This proves that $1 \to K(G) \to \Aut(G) \to \GL(d, \F_p) \to 1$ is exact.

Let $x_1, \dots, x_d$ be a minimal generating set for $G$.  Also let $L$ be a normal subgroup of $G$ lying in $G_2$ and let $u_1, \dots, u_d$ be any elements of $L$.  Since $G$ is relatively free, the map $\alpha : x_i \mapsto x_i u_i$ for each $i$ is an endomorphism of $G$ (it suffices to check that if a word $w$ in the $x_i$'s equals 1, then $w\alpha = 1$, but every tuple of elements of $G$ satisfies the same relations, so when $x_i$ is replaced by $x_i u_i$ in $w$, the new word also equals 1).  Furthermore, $\alpha$ acts trivially on $G/L$ and is an automorphism by Proposition~\ref{frataut}.  Conversely, any automorphism of $G$ that acts trivially on $G/L$ must act on each $x_i$ as multiplication by an element of $L$.  Thus the number of automorphisms of $G$ that act trivially on $G/L$ is $|L|^d$.  Taking $L = G_2$ gives $|K(G)| = p^{d(g-d)}$.

Next, we claim that any group $H$ in $V$ that is minimally generated by $d$ elements is isomorphic to $G/L$ for some normal subgroup $L$ of $G$ lying in $G_2$.  Evidently $H$ is isomorphic to $G/L$ for some normal subgroup $L$ of $G$; it suffices to show that if $L \not\le G_2$, then $G/L$ will be generated by fewer than $d$ elements.  Choose $x_1 \in L \setminus G_2$.  Extend $\{x_1\}$ to a generating set $\{x_1, \dots, x_d\}$ of $G$.  Then $G/L$ is generated by the images of $\{x_2, \dots, x_d\}$.

Suppose $M$ is a normal subgroup of $G$ in the same $\Aut(G)$-orbit as $L$.  Clearly $G/M \cong G/L$, so the map $L \mapsto G/L$ is well-defined on $\Aut(G)$-orbits of normal subgroups of $G$ lying in $G_2$.  To show that this is a bijection, we must show that if $M$ is a normal subgroup of $G$ lying in $G_2$ with $G/M \cong G/L$, then $M$ is in the same $\Aut(G)$-orbit as $L$.  Let $\beta : G/L \to G/M$ be an isomorphism.  By~\cite[Theorem 44.21]{neu}, $G$ is projective, as in \cite[Definition 44.11]{neu}; as the quotient map from $G$ to $G/M$ is surjective, this says that there exists an endomorphism $\gamma : G \to G$ so that the diagram in Figure 1 commutes.  Then $\gamma$ induces $\beta$, and $\beta$ induces an automorphism on the Frattini quotient of $G$ (since the Frattini quotients of $G/L$ and $G/M$ are isomorphic to the Frattini quotient of $G$).  It follows from Proposition~\ref{frataut} that $\gamma$ is an automorphism of $G$.  From Figure 1, it is also clear that $L\gamma \le M$.  Thus $L \gamma = M$, and $L$ and $M$ are in the same $\Aut(G)$-orbit.

\[
\begin{array}{c}
\xymatrix{
G \ar@{.>}[d]_{\gamma} \ar[r] & G/L \ar[d]^{\beta} \\
G \ar[r] & G/M
} \\
\textrm{Figure 1}
\end{array}
\]

If we take $L = M$, we find that any automorphism of $H = G/L$ is induced by an automorphism of $G$, so that $\Aut(H) \cong N_{\Aut(G)}(L)/B(L)$, where $B(L)$ is the subgroup of $N_{\Aut(G)}(L)$ that acts trivially on $H$.  By the earlier argument in this proof, $|B(L)| = |L|^d$.
\end{proof}

\begin{thm} \label{isoenum2}
Suppose that $G$ is the relatively free group on $d$ generators in a variety of groups $V$ and suppose that $G$ has lower $p$-length $n$.  The map $L \mapsto G/L$ defines a bijection between $\GL(d,\F_p)$-orbits on normal subgroups $L$ of $G$ lying in $G_n$ and groups $H$ in $V$ that are minimally generated by $d$ elements and satisfy $H/H_n \cong G/G_n$.  If $H = G/L$, then
\[
1 \to K(G)/B(L) \to \Aut(H) \to N_{\GL(d,\F_p)}(L) \to 1
\]
is exact, where $B(L)$ is the subgroup of $N_{\Aut(G)}(L)$ that acts trivially on $H$.  Moreover, $K(H)$ is the image of $K(G)/B(L)$ in $\Aut(H)$.
\end{thm}

\begin{proof}
$H/H_n \cong G/G_nL$ is isomorphic to $G/G_n$ if and only if $L \le G_n$.  Furthermore, $K(G)$ acts trivially on $G_n \cong G_n/G_{n+1}$ as noted in Subsection~\ref{autgroup}, so the $\Aut(G)$-orbits of normal subgroups of $G$ lying in $G_n$ are just the $\GL(d,\F_p)$-orbits.  This proves the bijection.

Since $K(G)$ fixes $L$, it also follows that 
\[
1 \to K(G) \to N_{\Aut(G)}(L) \to N_{\GL(d,
\F_p)}(L) \to 1
\]
is exact.  Combined with the second exact sequence in Theorem~\ref{isoenum1}, we find that
\[
1 \to K(G)/B(L) \to \Aut(H) \to N_{\GL(d,\F_p)}(L) \to 1
\]
is exact.  Every automorphism in $K(G)$ induces an automorphism in $K(H)$ since $K(G)$ fixes $L$ and $G/G_2 \cong H/H_2$.  Conversely, every automorphism in $K(H)$ is induced by an automorphism in $K(G)$.  The kernel of the map from $K(G)$ to $K(H)$ is $B(L)$, so $K(H)$ is the image of $K(G)/B(L)$.
\end{proof}

We can now prove Theorem~\ref{bijthm}, restated here for convenience.

\begin{thmx}[Theorem~\ref{bijthm}]
Fix a prime $p$ and integers $d, n \ge 2$.  Let $F$ be the free group on $d$ generators and define the following sets:
\begin{eqnarray*}
\mathcal{A}_{d,n} &=& \{\textrm{normal subgroups of $F/F_{n+1}$ lying in $F_2/F_{n+1}$}\} \\
\mathcal{B}_{d,n} &=& \{\textrm{normal subgroups of $F/F_{n+1}$ lying in $F_2/F_{n+1}$} \\
&& \qquad \textrm{and not containing $F_n/F_{n+1}$}\} \\
\mathcal{C}_{d,n} &=& \{\textrm{normal subgroups of $F/F_{n+1}$ lying in $F_n/F_{n+1}$}\} \\
\mathcal{D}_{d,n} &=& \{\textrm{normal subgroups of $F/F_{n+1}$ contained in the} \\
&& \qquad \textrm{regular $\GL(d,\F_p)$-orbits in $\mathfrak{C}_{d,n}$}\} \\
\\
\mathfrak{A}_{d,n} &=& \{\textrm{$\Aut(F/F_{n+1})$-orbits in $\mathcal{A}_{d,n}$}\} \\
\mathfrak{B}_{d,n} &=& \{\textrm{$\Aut(F/F_{n+1})$-orbits in $\mathcal{B}_{d,n}$}\} \\
\mathfrak{C}_{d,n} &=& \{\textrm{$\Aut(F/F_{n+1})$-orbits in $\mathcal{C}_{d,n}$}\} = \{\textrm{$\GL(d,\F_p)$-orbits in $\mathcal{C}_{d,n}$}\} \\
\mathfrak{D}_{d,n} &=& \{\textrm{regular $\GL(d,\F_p)$-orbits in $\mathcal{C}_{d,n}$}\}.
\end{eqnarray*}
Then there is a well-defined map $\pi_{d,n}: \mathfrak{A}_{d,n} \to \{\textrm{finite $p$-groups}\}$ given by $L/F_{n+1} \mapsto F/L$, where $L/F_{n+1} \in \mathcal{A}_{d,n}$.  Furthermore $\pi_{d,n}$ induces bijections
\begin{eqnarray*}
\mathfrak{A}_{d,n} &\leftrightarrow& \{\textrm{$p$-groups of lower $p$-length at most $n$} \\
&& \qquad \textrm{and minimally generated by $d$ elements} \} \\
\mathfrak{B}_{d,n} &\leftrightarrow& \{\textrm{$p$-groups of lower $p$-length $n$} \\
&& \qquad \textrm{and minimally generated by $d$ elements} \} \\
\mathfrak{D}_{d,n} &\leftrightarrow& \{\textrm{subgroups $H$ in $\pi_{d,n}(\mathfrak{C}_{d,n})$ with $A(H) = 1$}\}.
\end{eqnarray*}
\end{thmx}

\begin{proof}
Take $V$ to be the variety of $p$-groups of lower $p$-length at most $n$.  Then $F/F_{n+1}$ is the relatively free group on $d$ generators in $V$.  The $\Aut(F/F_{n+1})$- and $\GL(d,\F_p)$-orbits in $\mathcal{C}_{d,n}$ are the same because of the first exact sequence in Theorem~\ref{isoenum1} and the fact that $K(F/F_{n+1})$ acts trivially on $F_n/F_{n+1}$ as in Subsection~\ref{autgroup}.

The map $\pi_{d,n}$ is well-defined and defines bijections for $\mathfrak{A}_{d,n}$ and $\mathfrak{B}_{d,n}$ by Theorem~\ref{isoenum1}.  A normal subgroup $L$ of $F/F_{n+1}$ lying in $F_n/F_{n+1}$ is in a regular $\GL(d,\F_p)$-orbit if $N_{\GL(d,\F_p)}(L) = 1$.  By Theorem~\ref{isoenum2}, $L$ is in a regular orbit if and only if $A(H) = 1$.  Thus the bijection for $\mathfrak{D}_{d,n}$ is proved.
\end{proof}

Note, by the way, that since $F_n/F_{n+1}$ is elementary abelian and central in $F/F_{n+1}$, the set $\mathcal{C}_{d,n}$ is just the set of subspaces of the vector space $F_n/F_{n+1}$.

\section{The Lower $p$-Series of a Free Group}
\label{freegpsec}

Let $F$ be the free group on $d$ generators $y_1, y_2, \dots, y_d$.  To prepare for Sections~\ref{est1sec} and~\ref{est2sec}, we need to analyze the $\F_p \GL(d, \F_p)$-module structure of $F_n/F_{n+1}$ along with power and commutator maps from $F_n/F_{n+1}$ to $F_{n+1}/F_{n+2}$.  Our main tool will be the connection between the lower $p$-series of $F$ and the free Lie algebra described in Theorem~\ref{structurethm}.  The results of Theorem~\ref{structurethm} appear several times in the literature with varying degrees of correctness and detail. Our presentation follows Bryant and Kov\'{a}cs~\cite{bk}, while the most complete proofs may be inferred from Huppert and Blackburn~\cite[Chapter VIII]{hb}.  Information about the free Lie algebra can be found in Garsia~\cite{gar} and Reutenauer~\cite{reu}.

Let $K$ be any field and let $A = \{x_1, \dots, x_d\}$ be an alphabet on $d$ letters.  Write $A^{\ast}$ for the collection of all $A$-words and $A^n$ for the collection of all $A$-words of length $n$.  Let $K[A^{\ast}]$ denote the free associative $K$-algebra on the generators $x_1, x_2, \dots, x_d$; equivalently, $K[A^{\ast}]$ is the non-commutative algebra of polynomials
\[
f = \sum_{w \in A^{\ast}}{f_w w}
\]
with coefficients $f_w \in K$.  The algebra $K[A^{\ast}]$ is graded by degree; let $K[A^n]$ denote the homogeneous component of degree $n$.  Also, $K[A^{\ast}]$ is a Lie algebra under the Lie bracket $[f,g] = fg-gf$.  Let $K[\Lambda^{\ast}]$ denote the Lie subalgebra of $K[A^{\ast}]$ generated by $x_1, \dots, x_d$ and the Lie bracket.  Then $K[\Lambda^{\ast}]$ is the \emph{free Lie algebra} over $K$ on $x_1, \dots, x_d$.  It is also graded by degree; let $K[\Lambda^n]$ be the homogeneous component of $K[\Lambda^{\ast}]$ of degree $n$.

It will be convenient to specify a basis of $K[\Lambda^n]$.  Lexicographically order the set $A^{\ast}$, where $x_1 < x_2 < \cdots < x_d$.  A word $w$ is a \emph{Lyndon word} if it is smaller than all of its proper non-trivial tails.  Let $L$ be the set of Lyndon words, and let $L_n$ be the set of Lyndon words of length $n$.  Inductively define the \emph{right standard bracketing} $b[w]$ of $w \in L$ by
\[
b[w] = w
\]
if $w \in A$ and otherwise by
\[
b[w] = \left[b\left[w_1\right], b\left[w_2\right]\right],
\]
where $w = w_1w_2$ and $w_2$ is the longest proper tail of $w$ that is a Lyndon word.

\begin{thm}[{Reutenauer~\cite[Proof of Theorem 5.1]{reu}}]
\label{basisthm}
If $w \in L$, then
\[
b[w] = w + \sum_{w < v}{f_v v}
\]
for some $f_v \in K$.  The set $\{b[w] : w \in L_n\}$ forms a basis for $K[\Lambda^n]$.
\end{thm}

The results in this section require many maps; in an attempt to clarify matters, we will define all the maps now, using suggestive names, and postpone stating their properties until necessary.

\renewcommand{\arraystretch}{1.1}
\begin{defn}
Fix a prime $p$.  Fix integers $n \ge 1$, $d \ge 2$, and $1 \le j \le d$.  Let $f_i \in F_i$ for each $i \ge 1$.
\begin{itemize}
\item $\pow_n: F_n/F_{n+1} \to F_{n+1}/F_{n+2}$ 

(\emph{a power map on $F$})

\nopagebreak
\quad
$
\begin{array}{r@{\;:\;}c@{\quad\mapsto\quad}l}
\pow_n & f_n F_{n+1} & f_n^pF_{n+2}
\end{array}
$
\item $\Fcom_{j,n}: F_n/F_{n+1} \to F_{n+1}/F_{n+2}$ 

(\emph{a commutator map on $F$})

\nopagebreak
\quad
$
\begin{array}{r@{\;:\;}c@{\quad\mapsto\quad}l}
\Fcom_{j,n} & f_n F_{n+1} & [f_n, y_j]F_{n+2}
\end{array}
$
\item $\emb_n : F_n \to \F_p[A^{\ast}]$

(\emph{an embedding of $F_n$ into $\F_p[A^{\ast}]$})

\nopagebreak
\quad 
$
\begin{array}{r@{\;:\;}c@{\quad\mapsto\quad}l}
\emb_1 & y_j & x_j \\
\emb_n & f_{n-1}^p & \left\{ \begin{array}{c@{\;:\;}l}
		\emb_1(f_1) + \emb_1(f_1)^2 & \textrm{$n = 2$ and $p = 2$} \\
		\emb_{n-1}(f_{n-1}) & \textrm{otherwise}
	\end{array} \right.  \\
\emb_n & \left[ f_{n-1}, f_1\right] & \left[\emb_{n-1}(f_{n-1}), \emb_1(f_1) \right] \\
\emb_n & f_{n+1} & 0
\end{array}
$
\item $\qemb_n : F_n/F_{n+1} \to \F_p[A^{\ast}]$

(\emph{an embedding of the quotient $F_n/F_{n+1}$ into $\F_p[A^{\ast}]$})

\nopagebreak
\quad 
$
\begin{array}{r}
\textrm{$\qemb_n$ is induced by $\emb_n$}
\end{array}
$
\item $\com : \{\textrm{subspaces of $\F_p[A^{\ast}]$}\} \to \{\textrm{subspaces of $\F_p[A^{\ast}]$}\}$

(\emph{a commutator map on $\F_p[A^{\ast}]$})

\nopagebreak
\quad
$
\begin{array}{r@{\;:\;}c@{\quad\mapsto\quad}l}
\com & W & [W, \F_p[\Lambda^1]]
\end{array}
$
\item $\com_j : \F_p[A^{\ast}] \to \F_p[A^{\ast}]$

(\emph{a commutator map on $\F_p[A^{\ast}]$})

\nopagebreak
\quad
$
\begin{array}{r@{\;:\;}c@{\quad\mapsto\quad}l}
\com_j & f & [f, x_j]
\end{array}
$
\item $\com_{j,n} : \F_p[A^n] \to \F_p[A^{n+1}]$

(\emph{a commutator map on $\F_p[A^n]$})

\nopagebreak
\quad
$
\begin{array}{r}
\textrm{$\com_{j,n}$ is induced by $\com_j$.}
\end{array}
$
\item $\proj_n : \F_p[A^{\ast}] \to \F_p[A^n]$

(\emph{the projection map onto $\F_p[A^n]$})
\end{itemize}
\end{defn}
\renewcommand{\arraystretch}{1}

\begin{thm}
\label{structurethm}
The map $\emb_n$ is a well-defined homomorphism.  The map $\qemb_n$ is an $\F_p\GL(d, \F_p)$-module embedding of $F_n/F_{n+1}$ into $\F_p[A^{\ast}]$.  If $p$ is odd, the image of $\qemb_n$ is $\F_p[\Lambda^1] \oplus \cdots \oplus \F_p[\Lambda^n]$, and hence
\[
F_n/F_{n+1} \cong \F_p[\Lambda^1] \oplus \cdots \oplus \F_p[\Lambda^n]
\]
as $\F_p\GL(d, \F_p)$-modules.

If $p = 2$, the image of $\qemb_1$ is $\F_2[\Lambda^1]$.  The image $E$ of $\qemb_2$ satisfies
\[
E + \F_2[A^2] = \F_2[A^1] \oplus \F_2[A^2] \quad \textrm{and} \quad E \cap \F_2[A^2] = \F_2[\Lambda^2],
\]
so $E$ is an extension of $\F_2[\Lambda^2]$ by $\F_2[\Lambda^1]$.  For $n \ge 3$, the image of $\qemb_n$ is $E \oplus \F_2[\Lambda^3] \oplus \cdots \oplus \F_2[\Lambda^n]$, and hence
\[
F_n/F_{n+1} \cong E \oplus \F_2[\Lambda^3] \oplus \cdots \oplus \F_2[\Lambda^n].
\]
\end{thm}

Note that as a $\F_p\GL(d, \F_p)$-module, $\F_p[\Lambda^n] \cong V \wedge V \wedge \cdots \wedge V$, the $n$-fold wedge product where $V$ is the natural $\F_p\GL(d, \F_p)$-module.

\begin{cor}
\label{maps}
Unless $p = 2$ and $n = 1$, the diagram on the left commutes and $\pow_n$ is an injective homomorphism.  The diagram on the right commutes and $\Fcom_{j,n}$ is a homomorphism.
\[
\begin{array}{c}
\xymatrix{
F_n/F_{n+1} \ar[dr]_{\pow_n} \ar[rr]^{\qemb_n} && \ar[dl]^{\qemb_{n+1}} \F_p[A^{\ast}] \\
& F_{n+1}/F_{n+2}
}
\end{array} \qquad 
\begin{array}{c}
\xymatrix{
F_n/F_{n+1} \ar[d]_{\Fcom_{j,n}} \ar[r]^{\qemb_n} & \F_p[A^{\ast}] \ar[d]^{\com_j} \\
F_{n+1}/F_{n+2} & \ar[l]_{\qemb_{n+1}} \F_p[A^n] 
}
\end{array}
\]
\end{cor}

The dimension of $K[\Lambda^i]$ is given by Witt's formula:
\[
\dim(K[\Lambda^i]) = \frac{1}{i} \sum_{j | i}{\mu(i/j) \cdot d^j},
\]
where $\mu$ is the M\"{o}bius function (see~\cite[Appendix 0.4.2]{reu}).  Thus Theorem~\ref{structurethm} tells us the rank of $F_n/F_{n+1}$.

\begin{cor}
\label{dimcor}
The rank of $F_n/F_{n+1}$ is
\[
\sum_{i=1}^n{\frac{1}{i} \sum_{j | i}{\mu(i/j) \cdot d^j}}.
\]
\end{cor}

The remainder of this section is devoted to proving the following theorem and corollary.  Corollary~\ref{expcor} will allow us to count normal subgroups of $F/F_{n+1}$ when combined with Theorem~\ref{numnorms}.  

\begin{thm} 
\label{expthm}
Fix a prime $p$ and integers $d \ge 3$ and $n \ge 2$.  Suppose that $U$ is a normal subgroup of $F$ lying in $F_2$.  Let
\begin{eqnarray*}
Q &=& (U \cap F_n)F_{n+1}/F_{n+1} \\
R &=& (U_2 \cap F_{n+1})F_{n+2}/F_{n+2} \\
S &=& (U^p[U,F] \cap F_{n+1})F_{n+2}/F_{n+2}.
\end{eqnarray*}
Then $\rank(R) \ge \rank(Q)$ and $\rank(S) \ge (3/2) \; \rank(Q)$.
\end{thm}

The third isomorphism theorem lets us replace $F$ by $F/F_n$, giving the following corollary.

\begin{cor}
\label{expcor}
Fix a prime $p$ and integers $d \ge 3$, $n \ge 3$, and $2 \le i < n$.  Let $G = F/F_{n+1}$.  Suppose that $U$ is a normal subgroup of $G$ lying in $G_2$.  Let
\begin{eqnarray*}
Q &=& (U \cap G_i)G_{i+1}/G_{i+1} \\
R &=& (U_2 \cap G_{i+1})G_{i+2}/G_{i+2} \\
S &=& (U^p[U,F] \cap G_{i+1})G_{i+2}/G_{i+2}.
\end{eqnarray*}
Then $\rank(R) \ge \rank(Q)$ and $\rank(S) \ge (3/2) \; \rank(Q)$.
\end{cor}

To prove Theorem~\ref{expthm}, we will build up to an analogous result for the free Lie algebra on $d$ generators (Lemma~\ref{dimlem}) and then apply Theorem~\ref{structurethm}.

\begin{lem}
\label{injlem}
The following diagram commutes:
\[
\begin{array}{c}
\xymatrix{
\F_p[A^{\ast}] \ar[d]_{\proj_n} \ar[r]^{\com_j} & \F_p[A^{\ast}] \ar[d]^{\proj_{n+1}} \\
\F_p[A^n] \ar[r]_{\com_{j,n}} & \F_p[A^{n+1}]
}
\end{array}
\]
If $n = 1$, then the kernel of $\com_{j,n}$ is spanned by $x_j$.  If $n > 1$, then $\com_{j,n}$ is injective.
\end{lem}

\begin{proof}
The only statements requiring proof are those about the kernel and injectivity of $\com_{j,n}$.  Without loss of generality, we may assume that $j = 1$.  Suppose that $w \in L_n$.  Unless $n = 1$ and $w = x_1$, we see that $x_1 w$ is smaller than $w$, and hence smaller than all of its proper non-trivial tails.  So $x_1 w \in L_{n+1}$.  Furthermore, $w$ is the longest tail of $x_1 w$ that is a Lyndon word, so $b[x_1 w] = -[b[w], x_1]$.  Thus the image of $b[w]$ under $\com_{1,n}$ is the negative of a basis element in $L_{n+1}$, unique for each $w$.  It follows that the kernel of $\com_{1,1}$ is generated by $x_1$ and $\com_{1,n}$ is injective for $n > 1$.
\end{proof}

\begin{lem}\label{explem}
Fix $d \ge 3$ and $n \ge 2$.  Suppose that $W$ is a subspace of $\F_p[\Lambda^n]$.  Then $\dim(\com(W)) \ge (3/2)\dim(W)$.
\end{lem}

\begin{proof}
Let $\F_p[\Lambda^{\ast}]_{ij}$ denote the free Lie algebra on two generators $x_i$ and $x_j$; there is a natural embedding of $\F_p[\Lambda^{\ast}]_{ij}$ into $\F_p[\Lambda^{\ast}]$.  Let $\F_p[\Lambda^n]_{ij}$ be the homogeneous component of degree $n$ in $\F_p[\Lambda^{\ast}]_{ij}$.

First, we claim that if $f$ and $g$ are distinct elements of $\F_p[\Lambda^n]$ and $[f,x_i] = [g,x_j]$, then in fact $f, g \in \F_p[\Lambda^n]_{ij}$.  We may assume that $i,j > 1$.  Suppose that $f \notin \F_p[\Lambda^n]_{ij}$.  Then writing 
\[
f = \sum_{w \in L_n}{f_w b[w]},
\]
there must be some word $w \in L_n$ where $f_w \neq 0$ and $w$ contains a letter other than $x_i$ and $x_j$.  We may assume that $w$ contains the letter $x_1$.  In that case, by Theorem~\ref{basisthm}, there is a word beginning with $x_1$ that appears in $f$ with non-zero coefficient.  Thus there is a word beginning with $x_1$ and ending with $x_i$ that appears in $[f, x_i]$ with non-zero coefficient.  No such word can appear in $[g,x_j]$, contradicting the fact that $[f,x_i] = [g,x_j]$.  Hence $f \in \F_p[\Lambda^n]_{ij}$ and similarly $g \in \F_p[\Lambda^n]_{ij}$.

Note that $\F_p[\Lambda^n]_{ij} \cap \F_p[\Lambda^n]_{kl} = 0$ if $\{i,j\} \neq \{k,l\}$ (the letters $x_i$ and $x_j$ appear in every element of $\F_p[\Lambda^n]_{ij}$ since $n > 1$).  Choose $i$ and $j$ so that $\dim(W \cap \F_p[\Lambda^n]_{ij})$ is as small as possible; in particular this intersection has dimension at most $(1/2) \dim(W)$.  Let $X$ be a complement to $W \cap \F_p[\Lambda^n]_{ij}$ in $W$.

Define a more restrictive commutator map on subspaces by $\com_{ij}: \bullet \mapsto [\bullet, \F_p[\Lambda^1]_{ij}]$.  Obviously $\com_{ij}(W) \subseteq \com(W)$.  Using Lemma~\ref{injlem} and the above claim,
\begin{eqnarray*}
\dim(\com_{ij}(W)) &=& \dim(\com_{ij}(W \cap \F_p[\Lambda^n]_{ij})) + \dim(\com_{ij}(X)) \\
&\ge& \dim(W \cap \F_p[\Lambda^n]_{ij}) + 2 \dim(X) \\
&\ge& (3/2) \dim{W}.
\end{eqnarray*}
\end{proof}

\begin{lem}\label{onelem}
Fix $d \ge 2$.  Suppose that $W$ be a subspace of $\F_p[\Lambda^1]$.  Then
$\dim(W + \com(W)) \ge (3/2)\dim(W)$.
\end{lem}

\begin{proof}
Recalling Lemma~\ref{injlem}, this is clear if $\dim(W) = 1$, and otherwise 
\[
\dim(\com_{1,1}(W)) \ge \dim(W)-1,
\]
implying the result since $W$ and $\com(W)$ are disjoint.
\end{proof}

\begin{lem}\label{elem}
Let $p = 2$.  Suppose that $W$ is a subspace of $E$, where $E$ is defined in Theorem~\ref{structurethm}.  Then
$\dim(W + \com(W)) \ge (3/2)\dim(W)$.
\end{lem}

\begin{proof}
Let $X = W \cap \F_2[\Lambda^2]$ and let $Y$ be a complement to $X$ in $W$.  Note that $\dim(Y) = \dim(\proj_1(Y))$.  By Lemma~\ref{onelem}, 
\[
\dim(\proj_1(Y) + \com(\proj_1(Y))) \ge (3/2)\dim(\proj_1(Y)).
\]
By the commutative diagram in Lemma~\ref{injlem}, it follows that $Y + \com(Y)$ contains a subspace of dimension at least $(3/2)\dim(Y)$ that has trivial intersection with $\F_2[\Lambda^3]$.  By Lemma~\ref{explem}, $\com(X) \le \F_2[\Lambda^3]$ contains a subspace of dimension at least $(3/2) \dim(X)$.  
Then
\begin{eqnarray*}
\dim(W + \com(W)) \ge (3/2)\dim(X) + (3/2)\dim(Y) = (3/2)\dim(W).
\end{eqnarray*}
\end{proof}

\begin{lem}\label{dimlem}
Fix $d \ge 3$.  Let $U_n = \F_p[\Lambda^1] \oplus \cdots \oplus \F_p[\Lambda^n]$ if $p$ is odd or $U_n = E \oplus \F_2[\Lambda^3] \oplus \cdots \oplus \F_2[\Lambda^n]$ if $p = 2$.  Suppose that $W$ is a subspace of $\F_p[A^{\ast}]$ contained in $U_n$.  Then $\dim(W + \com(W)) \ge (3/2) \dim(W)$.
\end{lem}

\begin{proof}
The proof will be by induction on $n$.  When $p$ is odd and $n = 1$, Lemma~\ref{onelem} gives the result.  When $p = 2$ and $n = 2$, Lemma~\ref{elem} gives the result.  So assume that $p$ is odd and $n > 1$ or that $p = 2$ and $n > 2$.  Assume the result holds for $n-1$.  Let $X = W \cap U_{n-1}$.  By the inductive hypothesis, 
\[
\dim(X + \com(X)) \ge (3/2)\dim(X).
\]
Furthermore, $X + \com(X) \le U_n$.  Let $Y$ be a complement to $X$ in $W$.  By the commutative diagram in Lemma~\ref{injlem}, $\com(\proj_n(Y)) = \proj_{n+1}(\com(Y))$.  By the definition of $X$ and $Y$, $\dim(\proj_n(Y)) = \dim(Y)$.  By Lemma~\ref{explem},
\[
\dim(\proj_{n+1}(\com(Y))) \ge (3/2) \dim(\proj_n(Y)).
\]
Thus $\com(Y)$ contains a subspace of dimension at least $(3/2) \dim(\proj_n(Y))$ that has trivial intersection with $U_n$.  Therefore
\[
\dim(W + \com(W)) \ge (3/2) \dim(X) + (3/2) \dim(Y) = (3/2) \dim(W).
\]
\end{proof}

\begin{proof}[Proof of Theorem~\ref{expthm}]
Replacing $U$ by $(U \cap F_n)F_{n+1}$ does not change $Q$, $R$, or $S$, so we may assume that $F_{n+1} \le U \le F_n$.  Recall that by Corollary~\ref{maps}, $\pow_n$ is injective.   Since $\pow_n(Q) = R$, it follows that $\rank(R) \ge \rank(Q)$.

Also by Corollary~\ref{maps},
\[
S = \qemb_{n+1}^{-1}(\qemb_n(U) + (\com \circ \qemb_n)(U)).
\]
Since $\qemb_n$ is injective, and 
\[
\dim(\qemb_n(U) + (\com \circ \qemb_n)(U)) \ge (3/2) \dim(\qemb_n(U))
\]
by Lemma~\ref{dimlem}, it follows that $\rank(S) \ge (3/2) \; \rank(Q)$.
\end{proof}

\section{Numerical Estimates} \label{est0sec}

The purpose of this section is to prove several estimates needed in Sections~\ref{est1sec} and~\ref{est2sec}.  Most of the estimates involve Gaussian coefficients, and so we will begin with the relevant definitions and bounds on the Gaussian coefficients obtained by Wilf~\cite{wil}.

The \emph{Gaussian coefficient} (also called the \emph{$q$-binomial coefficient})
\[
\sbinom{n}{k}_q = \frac{(q^n-1) \cdots (q^n-q^{k-1})}{(q^k-1) \cdots (q^k-q^{k-1})}
\]
is the number of $k$-dimensional subspaces of a vector space of dimension $n$ over $\F_q$.  We shall be concerned with estimates for $\sbinom{n}{k}_q$ and for the \emph{Galois number}
\[
\mathcal{G}_n(q) = \sum_{k=0}^n{\sbinom{n}{k}_q},
\]
which is the total number of subspaces of a vector space of dimension $n$ over $\F_q$.  (A survey of these numbers is given by Goldman and Rota~\cite{gr}.)  First we need a technical lemma.

\begin{lem}
\label{polybound}
Let 
\[
C(q) = \sum_{r=-\infty}^{\infty}{q^{-r^2}}.
\]
Let $f(x) = -ax^2 + bx + c$ with $a > 0$, let $|q| > 1$, and set $A(q) = \sum_r{q^{f(r)}}$, where the sum is over all integers $r$ with $t \le r \le u$.  Then $A(q) \le C(q^a) q^{f(y)}$ for some $y \in [t, u]$.
\end{lem}

\begin{proof}
Suppose the maximum of $f(x)$ in $[t, u]$ occurs at $x = y$.  The global maximum of $f(x)$ occurs at $x = b/2a$, so one of three cases holds: $b/2a \le y = t$, $u = y \le b/2a$, or $t \le y = b/2a = u$.  In each case, for all $r \in [t,u]$,
\begin{eqnarray*}
&&-a(r-y)^2 - f(r) + f(y) \\
&=& -a(r-y)^2 - (-ar^2+br+c) + (-ay^2+by+c) \\
&=& (2ay-b)(r-y) \\
&\ge& 0.
\end{eqnarray*}
Thus
\begin{eqnarray*}
A(q) &=& q^{f(y)} \sum_{t \le r \le u}{q^{f(r)-f(y)}} \\
&\le& q^{f(y)} \sum_{t \le r \le u}{q^{-a(r-y)^2}} \\
&\le& q^{f(y)} \sum_{r = -\infty}^{\infty}{q^{-a(r-y)^2}},
\end{eqnarray*}
and it suffices to show that
\[
g(y) = \sum_{r = -\infty}^{\infty}{s^{-(r-y)^2}} \le g(0),
\]
where $s = q^a$.  This is a consequence of Jacobi's functional equation for the theta function
\[
\theta_3(z, w) = \sum_{r = -\infty}^{\infty}{e^{r^2 \pi i w} e^{2 r i z}},
\]
where $|e^{\pi iw}| < 1$.  Section 21.51 of Whittaker and Watson~\cite{ww} gives the functional equation
\[
\theta_3(z, w) = \frac{1}{\sqrt{-iw}} e^{z^2/\pi i w} \theta_3(z/w, -1/w),
\]
where $\sqrt{e^{i \theta}}$ denotes $e^{i\theta/2}$ for $0 \le \theta \le 2\pi$.  Now 
\begin{eqnarray*}
g(y) &=& s^{-y^2} \sum_{r = -\infty}^{\infty}{s^{-r^2} e^{-2ri(iy \log{s})}} \\
&=& s^{-y^2} \theta_3(-iy \log{s}, w),
\end{eqnarray*}
where $s^{-1} = e^{\pi i w}$ so that $\pi i w = -\log{s}$.  Hence
\begin{eqnarray*}
g(y) &=& \frac{s^{-y^2} \sqrt{\pi}}{\sqrt{\log{s}}} e^{y^2 \log{s}} \theta_3(-\pi y, -1/w) \\
&=& \sqrt{\frac{\pi}{\log{s}}} \sum_{r = -\infty}^{r = \infty}{e^{-r^2 \pi^2 / \log{s}} e^{-2ir\pi y}} \\
&=& \sqrt{\frac{\pi}{\log{s}}} (1 + 2 \sum_{r=1}^{\infty}{e^{-r^2 \pi^2 / \log{s}} \cos{2r\pi y}}) \\
&\le& \sqrt{\frac{\pi}{\log{s}}} (1 + 2 \sum_{r=1}^{\infty}{e^{-r^2 \pi^2 / \log{s}}}) \\
&=& g(0).
\end{eqnarray*}
\end{proof}

To obtain bounds for Gaussian coefficients, let
\begin{eqnarray*}
D(q) &=& \prod_{j=1}^{\infty}{(1-q^{-j})^{-1}} \\
S_n(q) &=& \sum_{k=0}^n{q^{k(n-k)}} = q^{n^2/4} \sum_{k=0}^n{q^{-(k-n/2)^2}}.
\end{eqnarray*}

Note that both $C(q)$ and $D(q)$ decrease to 1 as $q \to \infty$.  If $q \ge 2$, then $C(q) \le C(2) < 9/4$ and $D(q) \le D(2) < 7/2$.  The following estimates on Gaussian coefficients and Galois numbers were either obtained by Wilf~\cite{wil} or follow from his work.

\begin{lem} \label{qests}
Fix $q \ge 2$.  Then
\begin{eqnarray}
\sbinom{n}{k}_q &\le& D(q) q^{k(n-k)} \label{coefbounds} \\
D(q) q^{n^2/4-1/4} \left(2 - \frac{9 q^{(1-n)/2}}{2} \right) &\le& \mathcal{G}_n(q) \label{gnbounds} \\
&\le& S_n(q) D(q) \nonumber \\
&\le& C(q) D(q) q^{n^2/4} \nonumber
\end{eqnarray}
\end{lem}

\begin{proof}
Equation~\ref{coefbounds} and $\mathcal{G}_n(q) \le S_n(q) D(q)$ are proved in~\cite{wil}.  The inequality $S_n(q) \le C(q) q^{n^2/4}$ follows from Lemma~\ref{polybound}, taking $f(x) = x(n-x) = -x^2+nx$ and noting that $x(n-x) \le n^2/4$ for all $x$.  This proves $\mathcal{G}_n(q) \le C(q) D(q) q^{n^2/4}$.

The lower bound for $\mathcal{G}_n(q)$ is slightly more complicated, but it is easy to see from~\cite{wil}, Lemma~\ref{polybound}, and the definition of $S_n(q)$ that
\begin{eqnarray*}
\mathcal{G}_n(q) &\ge& S_n(q) - \frac{2S_{n-1}(q)+2q^{-2n}}{q-1} \\
&\ge& 2 q^{n^2/4-1/4} - \frac{2 C(q) q^{(n-1)^2/4}}{q-1} \\
&\ge& q^{n^2/4-1/4} \left( 2 - \frac{2 C(q) q^{(1-n)/2}}{q-1} \right) \\
&\ge& q^{n^2/4-1/4} \left( 2 - \frac{9q^{(1-n)/2}}{2} \right),
\end{eqnarray*}
where the last inequality uses the fact that $2C(q)/(q-1) < 9/2$.
\end{proof}

Next we shall prove Lemma~\ref{gaussprods}, which will be needed in Section~\ref{est1sec} to bound products of Gaussian coefficients, and we will finish with Lemma~\ref{quadbound}, which will be used in Section~\ref{est2sec}. 

\begin{lem} \label{gaussprods}
Fix a prime $p$ and integers $n \ge 3$ and $d \ge 6$ or $n \ge 10$ and $d \ge 5$.  Let $F$ be the free group on $d$ generators, and let $d_n$ be the rank of $F_n/F_{n+1}$.  For $1 \le i \le n-1$ and $0 \le u_i \le d_i$, let
\[
A_i(u_i) = \sum{\prod_{j=i}^{n-1}{p^{-(u_{j+1}-d_{j+1})(u_{j+1}-u_j/2)}}},
\]
where the sum is over all integers $u_{i+1}, \dots, u_n$ such that
\begin{eqnarray*}
0 \le &u_j& \le d_j \qquad \textrm{for } i+1 \le j \le n-2 \\
1 \le &u_{n-1}& \le d_{n-1} \\
2 \le &u_n& \le d_n.
\end{eqnarray*} 
Then for $1 \le i \le n-2$, 
\[
A_i(u_i) \le C(p)^{n-i} p^{-15/16+d_n^2/4+d_{n-1}-d_n/4} p^{-u_i(d_{i+1}-1)/2}.
\]
\end{lem}

\begin{proof}
First note that
\[
A_{n-1}(u_{n-1}) = \sum_{u_n=2}^{d_n}{p^{-(u_n-d_n)(u_n-u_{n-1}/2)}}.
\]
As a function of $u_n$, the expression $-(u_n-d_n)(u_n-u_{n-1}/2)$ is at most $(d_n-u_{n-1}/2)^2/4$, so that
\[
A_{n-1}(u_{n-1}) \le C(p) p^{(d_n-u_{n-1}/2)^2/4}
\]
by Lemma~\ref{polybound}.

The proof of the theorem is by backward induction on $i$.  Note that
\[
A_i(u_i) = \sum_{u_{i+1}}{p^{-(u_{i+1}-d_{i+1})(u_{i+1}-u_i/2)} A_{i+1}(u_{i+1})}.
\]
When $i = n-2$, using our bound on $A_{n-1}(u_{n-1})$ gives
\begin{eqnarray*}
&& A_{n-2}(u_{n-2}) \\
&\le& C(p) p^{d_n^2/4} \sum_{u_{n-1} = 1}^{d_{n-1}}{p^{u_{n-1}^2/16-u_{n-1}d_n/4 + (d_{n-1}-u_{n-1})(u_{n-1}-u_{n-2}/2)}} \\
&=& C(p) p^{d_n^2/4} \sum_{u_{n-1} = 1}^{d_{n-1}}{p^{-15u_{n-1}^2/16+(-d_n/4+u_{n-2}/2+d_{n-1})u_{n-1}-d_{n-1}u_{n-2}/2}}
\end{eqnarray*}
As a function of $u_{n-1}$, the polynomial
\[
-15u_{n-1}^2/16+(-d_n/4+u_{n-2}/2+d_{n-1})u_{n-1}-d_{n-1}u_{n-2}/2
\]
is maximized at
\[
u_{n-1} = 8 (-d_n/4+u_{n-2}/2+d_{n-1})/15.
\]
Computations show that this is at most 1 when $n \ge 3$ and $d \ge 6$ or $n \ge 10$ and $d \ge 5$.  So as $u_{n-1}$ ranges from $1$ to $d_{n-1}$, the polynomial is maximized at $u_{n-1} = 1$.  By Lemma~\ref{polybound} and the fact that $C(p^{15/16}) \le C(p)$,
\begin{eqnarray*}
A_{n-2}(u_{n-2}) &\le& C(p)^2 p^{d_n^2/4-15/16-d_n/4+d_{n-1}} p^{(1-d_{n-1})u_{n-2}/2}.
\end{eqnarray*}
This proves the theorem for the base case $i = n-2$.  By induction, for $i \le n-3$,
\begin{eqnarray*}
A_{i}(u_i) &=& \sum_{u_{i+1} = 0}^{d_{i+1}}{p^{-(u_{i+1}-d_{i+1})(u_{i+1}-u_i/2)} A_{i+1}(u_{i+1})} \\
&\le& C(p)^{n-i-1} p^{-15/16+d_n^2/4+d_{n-1}-d_n/4} \\
&& \qquad \cdot \sum_{u_{i+1} = 0}^{d_{i+1}}{p^{-(u_{i+1}-d_{i+1})(u_{i+1}-u_i/2)-u_{i+1}(d_{i+2}-1)/2}}.
\end{eqnarray*}
As a function of $u_{i+1}$, the polynomial
\begin{eqnarray*}
&& -(u_{i+1}-d_{i+1})(u_{i+1}-u_i/2)-u_{i+1}(d_{i+2}-1)/2) \\
&=& -u_{i+1}^2+(d_{i+1}+u_i/2-(d_{i+2}-1)/2)u_{i+1} - d_{i+1}u_i/i
\end{eqnarray*}
is maximized at
\[
u_{i+1} = (d_{i+1} + u_i/i - (d_{i+2}-1)/2))/2.
\]
Computations show that this is at most 1/2 for $d \ge 3$ and $i \ge 1$.  So as $u_{i+1}$ ranges from $0$ to $d_{i+1}$, the polynomial is maximized at $u_{i+1} = 0$.  Thus
\[
A_i(u_i) \le C(p)^{n-i} p^{-15/16+d_n^2/4+d_{n-1}-d_n/4} p^{-(d_{i+1}-1)u_i/i}
\]
and the result is proved by induction.
\end{proof}

\begin{lem} \label{quadbound}
Suppose that $\alpha_1, \dots, \alpha_s$ are positive integers with $n = \alpha_1 + \cdots + \alpha_s$.  Then
\begin{equation} \label{ub1}
\alpha_1^2 + \cdots + \alpha_s^2 \le (n-s+1)^2 + (s-1),
\end{equation}
and this bound is achieved when $\alpha_1 = \alpha_2 = \cdots = \alpha_{s-1} = 1$.  Furthermore, if $n \ge \ep + 1$ and $s \ge 2$, then
\begin{equation} \label{ub2}
\alpha_1^2 + \cdots + \alpha_s^2 + \ep s \le (n-1)^2 + 1 + 2\ep.
\end{equation}
\end{lem}

\begin{proof}
For Equation~\ref{ub1}, we use a simple induction argument.  It is clearly true for $s = 1$.  Suppose it is true up through $s$; we will prove it for $s+1$.
\begin{eqnarray*}
\alpha_1^2 + \cdots + \alpha_s^2 + \alpha_{s+1}^2 &\le&
(n-\alpha_{s+1}-s+1)^2 + (s-1) + \alpha_{s+1}^2 \\
&\le& (n-s+1-\alpha_{s+1})^2 + \alpha_{s+1}^2 + (s-1) \\
&\le& (n-s+1-1)^2 + 1^2 + (s-1) \\
&=& (n-s)^2 + s,
\end{eqnarray*}
proving Equation~\ref{ub1}.  As for Equation~\ref{ub2},
\begin{eqnarray*}
\alpha_1^2 + \cdots + \alpha_s^2 + \ep s &\le& (n-s+1)^2 + (s-1) + \ep s \\
&=& ((n-1)-(s-2))^2 + s - 1 + \ep s \\
&=& (n-1)^2 - 2(n-1)(s-2) + (s-2)^2 + s - 1 + \ep s \\
&\le& (n-1)^2 - (\ep + s - 1)(s-2) + (s-2)^2 + s - 1 + \ep s \\
&=& (n-1)^2 + 1 + 2 \ep,
\end{eqnarray*}
where the first inequality follows from Equation~\ref{ub1} and the second inequality follows from the fact that since $n \ge \ep + 1$ and $n \ge s$, we know that $n \ge (\ep + s + 1)/2$.
\end{proof}

\section{From Subgroups in $F_2/F_{n+1}$ to Subgroups in $F_n/F_{n+1}$} \label{est1sec}

The goal of this section is to prove Theorem~\ref{est1thm}, essentially showing that most $\GL(d, \F_p)$-orbits of normal subgroups of $F/F_{n+1}$ contained in $F_2/F_{n+1}$ are $\GL(d, \F_p)$-orbits of normal subgroups of $F/F_{n+1}$ contained in $F_n/F_{n+1}$.  We will prove Theorem~\ref{est1thm} by estimating the number of normal subgroups of $F/F_{n+1}$ contained in $F_2/F_{n+1}$.  Theorem~\ref{numnorms} offers a refined estimate on the number of normal subgroups of an arbitrary finite $p$-group.  Our estimate depends on certain parameters which are difficult to work out in general, but have been calculated for $F/F_{n+1}$ in Corollary~\ref{expcor}.  This will give us the tools to prove Theorem~\ref{est1thm}.

Let $H$ be a finite $p$-group of lower $p$-length $n$.  Given a normal subgroup $U$ of $H$, note that by the second isomorphism theorem,
\[
(U \cap H_i)/(U \cap H_{i+1}) \cong (U \cap H_i)H_{i+1}/H_{i+1},
\]
and this quotient is elementary abelian.  Let
\[
S(H, \vec{u}) = \{ U \lhd H \; : \; \dim((U \cap H_i)H_{i+1}/H_{i+1}) = u_i \},
\]
where $\vec{u} = (u_1, \dots, u_n)$ and each integer $u_i$ satisfies 
\[
0 \le u_i \le h_i = \dim(H_i/H_{i+1}).
\]

\begin{thm}
\label{numnorms}
Suppose that for each $U \in S(H, \vec{u})$, 
\[
\dim((U_2 \cap H_i)H_{i+1}/H_{i+1}) \ge v_i
\]
and
\[
\dim((U^p[U, H] \cap H_i)H_{i+1}/H_{i+1}) \ge w_i.
\]
Then
\[
|S(H, \vec{u})| \le \sbinom{h_1}{u_1}_p \prod_{i=2}^n{\sbinom{h_i-w_i}{u_i-w_i}_p p^{(u_1 + \cdots + u_{i-1}- v_1 - \cdots - v_{i-1})(h_i-u_i)}}.
\]
\end{thm}

\begin{proof}
The proof proceeds by induction on $n$, the lower $p$-length of $H$.  If $n = 1$, then $H$ is elementary abelian of dimension $h_1$, so that $\vec{u} = (u_1)$ and $S(H, \vec{u}) = \sbinom{h_1}{u_1}_p$.

Now suppose that the result holds in $J = H/H_n$, a group which has lower $p$-length $n-1$.  Any normal subgroup $U$ of $H$ lying in $S(H, \vec{u})$ determines the subgroup $K = U \cap H_n$ of $H_n$ and the normal subgroup $L = UH_n/H_n$ of $J$.  The subgroup $K$ contains $U^p [U, H] \cap H_n$, by hypothesis $\dim(U^p [U, H] \cap H_n) \ge w_n$, and $\dim(K) = \dim(U \cap H_n) = u_n$.

For $1 \le i \le n-1$, since $J_i = H_i/H_n$,
\begin{eqnarray}
(L \cap J_i)J_{i+1}/J_{i+1} &=& (UH_n/H_n \cap H_i/H_n)(H_{i+1}/H_n)/(H_{i+1}/H_n) \nonumber \\
&\cong& (UH_n \cap H_i)H_{i+1}/H_{i+1} \label{congeqn} \\
&\cong& (U \cap H_i)H_{i+1}/H_{i+1}. \nonumber
\end{eqnarray}
Thus $L \in S(J, \vec{t} \;)$, where $\vec{t} = (u_1, \dots, u_{n-1})$.  Furthermore, if $M$ is the inverse image of $L$ in $H$, then
\[
M^p [M, H] = (UH_n)^p [UH_n, H] = U^p [U, H],
\]
since $H_{n+1} = H_n^p[H_n, H] = 1$.  Thus $L$ determines $U^p [U, H] \cap H_n$.

Given $L$, the subgroup $K$ is a subspace of $H_n$ of dimension $u_n$ containing $M^p [M, H] \cap H_n$, which has dimension at least $w_n$.  Let $w = \dim(M^p [M, H] \cap H_n)$.  Then there are
\[
\sbinom{h_n-w}{u_n-w}_p = \sbinom{h_n-w}{h_n-u_n}_p
\]
choices for $K$.  This Gaussian coefficient is a decreasing function of $w$, so there are at most
\[
\sbinom{h_n-w_n}{u_n-w_n}_p
\]
choices for $K$.  Hence the number of possible pairs $K$ and $L$ given by subgroups in $S(H, \vec{u})$ is at most
\[
|S(J, \vec{t} \;)| \cdot \sbinom{h_n - w_n}{u_n - w_n}_p.
\]
There is a bijection between subgroups $U \in S(H, \vec{u})$ that give $K$ and $L$ and complements to $H_n/K$ in $M/K$, given by $U \mapsto U/K$.  In the one direction, $U/K$ is a complement to $H_n/K$ since $U \cap H_n = K$ and $UH_n/K = M/K$.  In the other direction, a complement $U/K$ to $H_n/K$ satisfies $U \cap H_n = K$ and $UH_n/K = M/K$, so $U$ gives $K$ and $L$.

Recall that in general, if $G$ is a group with normal subgroup $N$, then the number of complements to $N$ in $G$ is either 0 or $|\Der(G/N, N)|$.  When $N$ is central, $\Der(G/N, N) = \Hom(G/N, N)$, and if the number of complements is 0, then $\Hom(G/N, N)$ is trivial (see Lubotzky and Segal~\cite[Lemma 1.3.1]{ls}).

Since $H_n/K$ is central in $M/K$ ($H_n \in Z(H)$), the number of complements to $H_n/K$ in $M/K$ is
\[
|\Hom(M/H_n, H_n/K)| = |\Hom(L, H_n/K)| = |\Hom(L/L_2, H_n/K)|.
\]
The dimension of $H_n/K = H_n/(H_n \cap U)$ is $h_n - u_n$.  Also,
\begin{eqnarray*}
&& \dim(L/L_2) \\
&=& \dim(L) - \dim(L_2) \\
&=& \sum_{i=1}^{n-1}{\dim((L \cap J_i)J_{i+1}/J_{i+1})} - \sum_{i = 1}^{n-1}{\dim((L_2 \cap J_i)J_{i+1}/J_{i+1})}.
\end{eqnarray*}
Note that $L_2 = U_2H_n/H_n$, and a similar calculation to Equation~\ref{congeqn} shows that
\[
(L_2 \cap J_i)J_{i+1}/J_{i+1} \cong (L_2 \cap H_i)H_{i+1}/H_{i+1},
\]
which by hypothesis has dimension at least $v_i$.  Thus
\[
\dim(L/L_2) \le u_1 + \cdots + u_{n-1} - (v_1 + \cdots + v_{n-1})
\]
and
\[
|\Hom(L/L_2, H_n/K)| \le p^{(h_n - u_n)(u_1 + \cdots + u_{n-1} - v_1 - \cdots - v_{n-1})}.
\]
Using the inductive hypothesis gives
\begin{eqnarray*}
S(H, \vec{u}) &\le& S(J, \vec{t} \; ) \cdot \sbinom{h_n - u_n}{u_n - w_n}_p \cdot p^{(h_n - u_n)(u_1 + \cdots + u_{n-1} - v_1 - \cdots - v_{n-1})} \\
&\le& \sbinom{h_1}{u_1}_p \prod_{i=2}^n{\sbinom{h_i-w_i}{u_i-w_i}_p p^{(u_1 + \cdots + u_{i-1}- v_1 - \cdots - v_{i-1})(H_i-u_i)}}.
\end{eqnarray*} 
\end{proof}

We can now prove Theorem~\ref{est1thm}, restated here for convenience.

\begin{thmx}[Theorem~\ref{est1thm}]
Fix a prime $p$ and integers $d$ and $n$ so that either $n \ge 3$ and $d \ge 6$ or $n \ge 10$ and $d \ge 5$.  Let $F$ be the free group on $d$ generators and let $d_n$ be the rank of $F_n/F_{n+1}$.  Then
\[
1 \le \frac{|\mathfrak{A}_{d,n}|}{|\mathfrak{C}_{d,n}|} \le 1 + C(p)^{n-1} D(p)^{n-2} p^{d_{n-1} - d_n/4 + d^2}.
\]
\end{thmx}

\begin{proof}
To prove this result, we need to apply the estimates of Lemmas~\ref{qests} and~\ref{gaussprods} to the upper bound for $S(H, \vec{u})$ obtained in Theorem~\ref{numnorms} in the case when $H = F/F_{n+1}$.  By Corollary~\ref{expcor}, we may choose $v_{i+1} = u_i$ and $w_{i+1} = (3/2)u_i$.  In particular, $w_{i+1} = 0$ if $u_i = 0$.  By Equation~\ref{coefbounds} of Lemma~\ref{qests}, we have
\[
\sbinom{n}{k}_p \le D(p) p^{k(n-k)}.
\]
Substituting in the bound obtained in Theorem~\ref{numnorms}, we find that, if $u_1 = 0$, then
\[
|S(H, \vec{u})| \le D(p)^{n-1} p^h,
\]
where
\begin{eqnarray*}
h &=& u_2(d_2 - u_2) + (u_3 - w_3)(d_3 - u_3) + \cdots + (u_n - w_n)(d_n - u_n) \\
&& \hspace{1.5in} + u_2(d_3-u_3) + \cdots + u_{n-1}(d_n - u_n) \\
&\le& -(u_2-d_2)u_2 - (u_3 - d_3)(u_3 - u_2/2)
- \cdots \\
&& \hspace{1.5in} -(u_n - d_n)(u_n - u_{n-1}/2).
\end{eqnarray*}
Hence
\[
|\mathcal{A}_{d,n}| \le |\mathcal{C}_{d, n}| + \sum_{\vec{u}}{D(p)^{n-1} p^h} = |\mathcal{C}_{d, n}| + D(p)^{n-1} \sum_{\vec{u}}{p^h},
\]
where the sum is taken over all $\vec{u}$ such that $U \in S(H,\vec{u})$ if and only if $U \le F_2/F_{n+1}$ and $U \not\le F_n/F_{n+1}$.  In terms of $\vec{u}$, this means that $u_{n-1} \ge 1$ and $u_1 = 0$.  Since $u_n \ge w_n > u_{n-1}$, we know that $u_n \ge 2$.  Then by Lemma~\ref{gaussprods}, we have
\[
\sum_{\vec{u}}{p^h} = A_1(0),
\]
and
\[
|\mathcal{A}_{d, n}| \le |\mathcal{C}_{d, n}| + D(p)^{n-1} C(p)^{n-1} p^y, 
\]
where 
\[
y = d_n^2/4-15/16-d_n/4+d_{n-1}.
\]
Hence, as $|\mathcal{C}_{d, n}| = \mathcal{G}_{d_n}(p)$, using  Lemma~\ref{qests} and the fact that $2 - 9p^{(1-d_n)/2}/2 > 1$,
\begin{eqnarray*}
|\mathcal{A}_{d,n}|/|\mathcal{C}_{d,n}| &\le& 1 + D(p)^{n-1} C(p)^{n-1} p^y / \mathcal{G}_{d_n}(p) \\
&\le& 1 + D(p)^{n-2} C(p)^{n-1} p^{d_{n-1} - d_n/4}.
\end{eqnarray*}
Now by Theorems~\ref{isoenum1} and~\ref{isoenum2}, $|\mathfrak{A}_{d,n}|$ and $|\mathfrak{C}_{d,n}|$ are the number of $\GL(d,\F_p)$-orbits on $\mathcal{A}_{d, n}$ and $\mathcal{C}_{d,n}$ respectively.  Hence
\[
0 \le |\mathfrak{A}_{d,n}| - |\mathfrak{C}_{d,n}| \le |\mathcal{A}_{d,n}| - |\mathcal{C}_{d,n}|,
\]
since $|\mathfrak{A}_{d,n}| - |\mathfrak{C}_{d,n}|$ is the number of $\GL(d,\F_p)$ orbits in $\mathcal{A}_{d, n} \setminus \mathcal{C}_{d,n}$.  Also $|\mathcal{C}_{d,n}| \le |\mathfrak{C}_{d,n}| \cdot |\GL(d,\F_p)|$, since $\mathcal{C}_{d,n}$ falls into $|\mathfrak{C}_{d,n}|$ orbits, each of size at most $|\GL(d,\F_p)|$.  Then
\begin{eqnarray*}
0 &\le& \frac{|\mathfrak{A}_{d,n}|}{|\mathfrak{C}_{d,n}|} - 1 \\
&=& \frac{|\mathcal{C}_{d,n}|}{|\mathfrak{C}_{d,n}|} \left( \frac{|\mathfrak{A}_{d,n}| - |\mathfrak{C}_{d,n}|}{|\mathcal{C}_{d,n}|} \right) \\
&\le& |\GL(d,\F_p)| \left( \frac{|\mathcal{A}_{d,n}| - |\mathcal{C}_{d,n}|}{|\mathcal{C}_{d,n}|} \right) \\
&\le& C(p)^{n-1} D(p)^{n-2} p^{d_{n-1} - d_n/4 + d^2}.
\end{eqnarray*}
Therefore
\[
1 \le \frac{|\mathfrak{A}_{d,n}|}{|\mathfrak{C}_{d,n}|} \le 1 + C(p)^{n-1} D(p)^{n-2} p^{d_{n-1} - d_n/4 + d^2}.
\]
\end{proof}

\section{Most Orbits on Subgroups of $F_n/F_{n+1}$ are Regular} \label{est2sec}

In this section we shall prove Theorem~\ref{est2thm}.  This depends on estimating $|\mathfrak{C}_{d,n}|$, the number of $\GL(d,\F_p)$-orbits on subspaces of $F_n/F_{n+1}$, via the Cauchy-Frobenius Lemma.  To do this, we obtain in Theorem~\ref{upperbound} an upper bound for the number of subspaces of $F_n/F_{n+1}$ fixed by an element of $\GL(d, \F_p)$, and refine this in Theorem~\ref{strongerbound} to obtain a stronger bound in the case $n = 2$.

Suppose $M$ is an $\F_p \GL(d, \F_p)$-module.  Let $g \in \GL(d, \F_p)$.  We want to count the number of subspaces of $M$ (viewed as an $\F_p$-vector space) fixed by $g$, which is the number of submodules of $M$ as a $\F_p \left< g \right>$-module.  We note that when $M$ is the natural $\F_p \GL(d, \F_p)$-module, Eick and O'Brien~\cite{eo} give an explicit formula for this number.  The following preliminaries are based on Macdonald~\cite[Chapter IV, Section 2]{mac}.  

Let $\Phi$ be the set of all polynomials in $\F_p[t]$ which are irreducible over $\F_p$ and let $P$ be the set of all partitions of non-negative integers.  Let $U$ be the set of all functions $\mu : \Phi \to P$ such that $m = \sum_{f \in \Phi}{\deg(f) |\mu(f)|}$, where $|\mu(f)|$ is the sum of the parts of the partition $\mu(f)$.  Then there is a one-to-one correspondence between $\F_p \left< g \right>$-modules $M$ of dimension $m$ and functions $\mu \in U$.  This correspondence is given by 
\begin{eqnarray*}
M &\cong& \bigoplus_{f \in \Phi}{\bigoplus_i{\frac{\F_p[t]}{(f)^{\mu_i(f)}}}},
\end{eqnarray*}
where $\mu_i(f)$ is the $i$-th part of $\mu(f)$, $(f)$ is the ideal of $\F_p[t]$ generated by $f$, and $g$ acts upon $\F_p[t]/(f)^s$ as multiplication by $t$.

Let
\[
M_f = \bigoplus_i{\frac{\F_p[t]}{(f)^{\mu_i(f)}}}.
\]
We call $\mu(f)$ the \emph{type} of $M_f$.  Any submodule $N$ of $M$ can be written $N = \oplus_{f \in \Phi}{N_f}$ with $N_f \subseteq M_f$ for each $f \in \Phi$.  That is, every submodule of $M$ is the direct sum of submodules of the summands $M_f$.  By Macdonald~\cite[Chapter II, 3.1]{mac} the type $\lambda$ of any $\F_p \left< g \right>$-submodule or quotient module of $M_f$ satisfies $\lambda \subseteq \mu(f)$.

For each $f \in \Phi$, let $\F_p[t]_f$ denote the localization of $\F_p[t]$ at the prime ideal $(f)$.  Then $\F_p[t]_f$ is a discrete valuation ring with residue field of order $q = p^{\deg(f)}$ and $M_f$ is a finite $\F_p[t]_f$-module of type $\mu(f)$.  

Both Theorems~\ref{upperbound} and \ref{strongerbound} depend on Theorem~\ref{submodules}, where we calculate the number of submodules of fixed type in a module of fixed type over a discrete valuation ring.  This generalizes the formula for the number of subgroups of a finite abelian $p$-group (see Birkhoff~\cite{bir}).

\begin{thm} \label{submodules}
Let $\mathfrak{a}$ be a discrete valuation ring with maximal ideal $\mathfrak{p}$ and let $\mathfrak{k} = \mathfrak{a}/\mathfrak{p}$ be the residue field of order $q$.  Let $\alpha = (\alpha_1, \alpha_2, \dots, \alpha_s)$ and $\beta = (\beta_1, \beta_2, \dots, \beta_r)$ be partitions with $\beta \subseteq \alpha$ and let $M$ be a finite $\mathfrak{a}$-module of type $\alpha'$.  Then the number of submodules of $M$ of type $\beta'$ is
\[
S(\alpha', \beta', q) = \prod_{i=1}^r{\sbinom{\alpha_i - \beta_{i+1}}{\beta_i - \beta_{i+1}}_q q^{\beta_{i+1} (\alpha_i - \beta_i)}}.
\]
\end{thm}

\begin{proof}
The proof is by induction on $\beta_1$.  If $\beta_1 = 0$, then $S(\alpha', \beta', q) = 1$ and the result holds.  Suppose $\beta_1 > 0$, and let the smallest part of $\beta'$ be $t$, so that either $\beta_1 = \cdots = \beta_t > \beta_{t+1}$ and $t < s$, or $\beta_1 = \cdots = \beta_s$ and $t = s$.  Write 
\[
\overline{\beta} = (\beta_1 - 1, \beta_2 - 1, \dots, \beta_t - 1, \beta_{t+1}, \dots).
\]
Let $N$ be any submodule of $M$ of type $\overline{\beta}'$, and le t$x$ be any element of $M$ with $\mathfrak{p}^t x = 0$, $\mathfrak{p}^{t-1} x \neq 0$, and $\mathfrak{p} x \cap N = 0$.  Then $\left< N, x \right>$ has type $\beta'$.  There are $S(\alpha', \overline{\beta}', q)$ choices for $N$, and for each $N$ it follows from~\cite[Chapter II, Equation 1.8]{mac} that the number of choices for $x$ is just
\begin{equation}
q^{\alpha_1 + \cdots + \alpha_t} ( 1 - q^{\beta_t - \alpha_t - 1}).
\end{equation}
On the other hand, fix a submodule $L$ of $M$ of type $\beta'$; we can count the number of choices of $N$ and $x$ so that $L = \left< N, x \right>$.  Here $N$ is a submodule of $L$ of type $\overline{\beta}'$ whose quotient has type $(t)$, and by~\cite[Chapter II, Equation 4.13]{mac}, the number of choices for $N$ is
\begin{eqnarray*}
\frac{1 - q^{\beta_{t+1} - \beta_t}}{1 - q^{-1}} \; q^{\sum_i{\binom{\beta_i}{2}} - \sum_i{\binom{\overline{\beta}_i}{2}}}
&=& \frac{1 - q^{\beta_{t+1} - \beta_t}}{1 - q^{-1}} \; q^{t (\beta_t - 1)}.
\end{eqnarray*}
Given $N$, it follows from~\cite[Chapter II, Equation 1.8]{mac} that there are 
\[
q^{\beta_1 + \cdots + \beta_t}(1-q^{-1})
\]
choices for $x$.  Thus any submodule $L$ of $M$ of type $\beta'$ arises as $\left< N, x \right>$ in 
\[
q^{\beta_1 + \cdots + \beta_t + t(\beta_t - 1)} (1 - q^{\beta_{t+1} - \beta_t})
\]
ways.  The total number of submodules $L$ of $M$ of type $\beta'$ is then
\begin{eqnarray}
S(\alpha', \beta', q) &=& \frac{S(\alpha', \overline{\beta}', q) q^{\alpha_1 + \cdots + \alpha_t} ( 1 - q^{\beta_t - \alpha_t - 1})}{
q^{\beta_1 + \cdots + \beta_t + t(\beta_t - 1)} (1 - q^{\beta_{t+1} - \beta_t})} \nonumber \\
&=& \frac{S(\alpha', \overline{\beta}', q) 
q^{\alpha_1 + \cdots + \alpha_t} ( 1 - q^{\beta_t - \alpha_t - 1})}{
q^{2t \beta_t - t} (1 - q^{\beta_{t+1} - \beta_t})}, \label{indexp}
\end{eqnarray}
where the second inequality uses $\beta_1 = \cdots = \beta_t$.  By induction, we know that
\begin{eqnarray*}
S(\alpha', \overline{\beta}', q) &=& \prod_{i=1}^r{\sbinom{\alpha_i - \overline{\beta}_{i+1}}{\overline{\beta}_i - \overline{\beta}_{i+1}}_q q^{\overline{\beta}_{i+1} (\alpha_i - \overline{\beta}_i)}} \\
&=& \prod_{i=1}^{t-1}{\sbinom{\alpha_i - \beta_{i+1} + 1}{\beta_i - \beta_{i+1}}_q q^{(\beta_{i+1} - 1)(\alpha_i - \beta_i + 1)}} \\
&& \qquad \cdot  \sbinom{\alpha_t - \beta_{t+1}}{\beta_t - \beta_{t+1} - 1}_q q^{\beta_{t+1} (\alpha_t - \beta_t + 1)} \\
&& \qquad \cdot \prod_{i=t+1}^r{\sbinom{\alpha_i - \beta_{i+1}}{\beta_i - \beta_{i+1}}_q q^{\beta_{i+1} (\alpha_i - \beta_i)}} \\
&=& \prod_{i=1}^r{\sbinom{\alpha_i - \beta_{i+1}}{\beta_i - \beta_{i+1}}_q q^{\beta_{i+1} (\alpha_i - \beta_i)}} \\
&& \qquad \cdot \prod_{i=1}^{t-1}{ \frac{q^{\alpha_i-\beta_{i+1}+1}-1}{q^{\alpha_i-\beta_i+1}-1} q^{\beta_{i+1} + \beta_i - \alpha_i - 1}} \cdot \frac{q^{\beta_t-\beta_{t+1}}-1}{q^{\alpha_t-\beta_t+1}-1} q^{\beta_{t+1}} \\
&=& \prod_{i=1}^r{\sbinom{\alpha_i - \beta_{i+1}}{\beta_i - \beta_{i+1}}_q q^{\beta_{i+1} (\alpha_i - \beta_i)}} \\
&& \qquad \cdot q^{2(t-1)\beta_t - \alpha_1 - \cdots - \alpha_{t-1} - (t-1)} \cdot \frac{q^{\beta_t-\beta_{t+1}}-1}{q^{\alpha_t-\beta_t+1}-1} q^{\beta_{t+1}} \\
&=& \prod_{i=1}^r{\sbinom{\alpha_i - \beta_{i+1}}{\beta_i - \beta_{i+1}}_q q^{\beta_{i+1} (\alpha_i - \beta_i)}} \cdot \frac{q^{2t \beta_t}}{q^{\alpha_1 + \cdots + \alpha_t + t}} \cdot \frac{1 - q^{\beta_{t+1} - \beta_t}}{1 - q^{\beta_t - \alpha_t - 1}}.
\end{eqnarray*}
Substituting this expression into Equation~\ref{indexp} gives the result.
\end{proof}

Using Theorem~\ref{submodules} and the techniques of Section~\ref{est0sec}, we can give an upper bound for the total number of submodules of a finite $\F_p \left< g \right>$-module $M$.  Note that every subspace of $M$ is a $\F_p \left< g \right>$-module if and only if $g$ acts as a scalar on $M$, that is, as multiplication by an element of $\F_p$.

\begin{thm} \label{upperbound}
Fix $d \ge 2$ and $g \in \GL(d, \F_p)$.  Suppose that $M$ is an $\F_p \left< g \right>$-module.  Let $m = \dim_{\F_p}(M)$ and let $S_M$ be the number of submodules of $M$.  Then either $g$ acts as a scalar on $M$ and $S_M = \mathcal{G}_m(p)$, or $g$ does not act as a scalar and
\[
\log_p{S_M} \le (m^2-2m+2)/4 + 2\ep,
\]
where $\ep = \log_p(C(p)D(p))$.
\end{thm}

\begin{proof}
Write $M = \oplus_{i=1}^k{M_i}$, where for each $i$, $M_i = M_{f_i}$ for some $f_i \in \Phi$ and $\dim_{\F_p}{M_i} = m_i$.  

\hspace{1in}

\noindent
\emph{Case 1: $k \ge 2$}.

Each submodule of $M$ is a direct sum of submodules of the summands $M_i$, so $S_M = \prod_{i=1}^k{S_{M_i}} \le \mathcal{G}_{m_1}(p) \mathcal{G}_{m-m_1}(p)$.  Then by Lemma~\ref{qests}, 
\[
S_M \le C(p)^2 D(p)^2 p^{m_1^2/4 + (m-m_1)^2/4} \le C(p)^2 D(p)^2 p^{(m^2-2m+2)/4},
\]
since $0 < m_1 < m$.

\hspace{1in}

\noindent
\emph{Case 2: $k = 1$}.

In this case, $M = M_f$ for some $f \in \Phi$.  Let $u = \deg(f)$ and $q = p^u$, and let $M$ have type $\alpha'$ as a $\F_p[t]_f$-module, where $\alpha = (\alpha_1, \dots, \alpha_s)$.  

\hspace{1in}

\noindent
\emph{Subcase 2.1: $\alpha$ has at least two parts}.

If $\beta = (\beta_1, \dots, \beta_r)$ and $\beta \subseteq \alpha$, then by Theorem~\ref{submodules} and Lemma~\ref{qests} Equation~\ref{coefbounds}, the number of submodules of $M$ of type $\beta'$ is
\begin{eqnarray*}
S(\alpha', \beta', q) &\le& \prod_{i=1}^r{D(q) q^{(\beta_i - \beta_{i+1})(\alpha_i - \beta_i) + \beta_{i+1}(\alpha_i - \beta_i)}} \\
&=& D(q)^r \prod_{i=1}^r{q^{\beta_i (\alpha_i - \beta_i)}}.
\end{eqnarray*}
Thus
\begin{eqnarray*}
S_M &=& \sum_{\beta' \subseteq \alpha'}{S(\alpha', \beta', q)} \\
&\le& D(q)^s \sum_{\beta' \subseteq \alpha'}{\prod_{i=1}^r{q^{\beta_i (\alpha_i - \beta_i)}}} \\
&\le& D(q)^s \prod_{i=1}^s{\sum_{b_i=0}^{\alpha_i}{q^{b_i(\alpha_i-b_i)}}} \\
&\le& D(q)^s C(q)^s \prod_{i=1}^s{q^{\alpha_i^2/4}},
\end{eqnarray*}
where the last inequality follows from Lemma~\ref{polybound}.  Now $D(q) \le D(p)$ and $C(q) \le C(p)$ so, remembering that $u(\alpha_1 + \cdots + \alpha_s) = m$ and using Lemma~\ref{quadbound},
\begin{eqnarray}
\log_p{S_M} &\le& u(\alpha_1^2 + \cdots + \alpha_s^2)/4 + s\ep \label{logbound} \\ 
&\le& (4 s \ep + (u\alpha_1)^2 + \cdots + (u\alpha_s)^2 + 4s \ep)/4\nonumber \\
&\le& ((m-1)^2 + 1 + 8 \ep)/4\nonumber \\
&\le& (m^2-2m+2)/4 + 2\ep, \nonumber
\end{eqnarray}
if $m \ge 4\ep + 1$.  For $m < 4\ep + 1$,   
\begin{eqnarray*}
\log_p{S_M} &\le& m^2/4 \\
&\le& (m^2-2m+2)/4 + (m-1)/2 \\
&\le& (m^2-2m+2)/4 + 2\ep.
\end{eqnarray*}

\hspace{1in}

\noindent
\emph{Subcase 2.2: $\alpha$ has one part}.

In this case, $\alpha_1 = m/u$.  If $u \ge 2$, then by Lemma~\ref{qests} Equation~\ref{coefbounds},
\begin{eqnarray*}
S_M &=& \sum_{0 \le \beta_1 \le \alpha_1}{\sbinom{\alpha_1}{\beta_1}_q} \\
&\le& C(q) D(q) q^{m^2/4u^2} \\
&\le& C(p)^2 D(p)^2 p^{m^2/4u} \\
&\le& C(p)^2 D(p)^2 p^{(m^2-2m+2)/4},
\end{eqnarray*}
since $u \ge 2$.  On the other hand, if $u = 1$, then $f = t-c$ for some $c \in \F_p$ and $M \cong \oplus^m \{\F_p[t]/(f)\}$ so that $g$ acts as the scalar $c$ on $M$ and $S_M = \mathcal{G}_m(p)$.
\end{proof}

The next theorem strengthens this result when the module structure is known more precisely and will be needed to deal with groups of lower $p$-length 2.

\begin{thm} \label{strongerbound}
Fix $d \ge 2$ and $g \in \GL(d, \F_p)$ with $g \neq 1$.  Suppose that $V$ is an $\F_p \left< g \right>$-module on which $g$ acts non-trivially and that $M$ is an $\F_p \left< g \right>$-module extension of $V \wedge V$ by $V$.  Let $v = \dim_{\F_p}(V)$, let $m = \dim_{\F_p}(M) = v(v+1)/2$, and let $S_M$ be the number of submodules of $M$.  Then 
\[
\log_p{S_M} \le (m-4)^2/4 + C,
\] 
where $\ep = \log_p{(C(p)D(p))}$ and
\[
C = \left\{
	\begin{array}{c@{\quad:\quad}l}
		\ep + 2m-4 & m \le 45 \\
		5\ep + 4 & \textrm{otherwise}.
	\end{array}
		\right.
\]
\end{thm}

\begin{proof}
First, if $v \le 9$, then $m \le 45$.  In this case,
\begin{eqnarray*}
S_M &\le& \mathcal{G}_m(p) \\
&\le& C(p) D(p) p^{m^2/4} \\
&=& C(p) D(p) p^{(m-4)^2/4 + 2m-4},
\end{eqnarray*}
proving the result.  So we may assume that $v \ge 10$.

Write $M = \oplus_{i=1}^k{M_i}$, where for each $i$, $M_i = M_{f_i}$ for some $f_i \in \Phi$ and $\dim_{\F_p}{M_i} = m_i$; we may assume that $m_1 \ge m_2 \ge \cdots \ge m_k$.  Note that $m_1 + \cdots + m_k = m$.  Then $V = \oplus_{i=1}^m{M_i \pi}$ where $\pi$ is the projection from $M$ onto $V$. 

Fix $0 < t < k$ and set $W = M_1 \oplus \cdots \oplus M_t$.  Also let $w = \dim{W} = m_1 + \cdots + m_t$.  Then $S_M \le \mathcal{G}_w(p) \mathcal{G}_{M-w}(p)$ since any submodule of $M$ is a direct sum of submodules of the summands $M_i$.  By Lemma~\ref{qests},
\[
S_M \le C(p)^2 D(p)^2 p^{w^2/4 + (M-w)^2/4}.
\]
When $4 \le w \le M-4$, it follows that
\begin{eqnarray*}
S_M &\le& C(p)^2 D(p)^2 p^{4 + (M-4)^2/4} \textrm{ and} \\
\log_p{S_M} &\le& (M-4)^2/4 + 2\ep + 4,
\end{eqnarray*}
proving the result.  If we cannot choose $t$ so that $4 \le w \le m-4$, then since $m > 9$ implies that $m_1 \not\le 3$, it must be that $m_1 \ge m-3$ and $k \le 4$.  Write $Y = M_2 \oplus \cdots \oplus M_k$; then $y = \dim{Y} \le 3$.  (It is possible that $Y$ is the zero module and that $y = 0$.)  At this point we need to prove a technical claim which we will use twice.

\hspace{1in}

\textbf{Claim:}  Suppose that $V$ is the direct sum of $\F_p \left< g \right>$-modules $A$ and $B$ of dimensions $a \ge 4$ and $v-a$ over $\F_p$, and suppose that $A \subset M_1 \pi$.  If $g$ acts as a scalar $c$ on $A$, then $c = 1$ and $A \otimes B$ is the direct sum of $a$ copies of $B$.

\hspace{1in}

\emph{Proof of claim:}  If $V = A \oplus B$, then $V \wedge V \cong (A \wedge A) \oplus (B \wedge B) \oplus (A \otimes B)$.  If $g$ acts as a scalar $c$ on $A$, then $A \cong \oplus \{ \F_p[t]/(t-c) \}^a$ and $M_1 = M_{f_1}$ with $f_1 = t-c$.  In this case $g$ acts as the scalar $c^2$ on $A \wedge A$, so $A \wedge A \cong \{ \F_p[t]/(t-c^2) \}^{a(a-1)/2}$.  If $c \neq 1$, then $A \wedge A \not\subseteq M_1$ and hence $A \wedge A \subseteq Y$.  But then $a(a-1)/2 = \dim(A \wedge A) \le \dim{Y} \le 3$, which is impossible.  Therefore $c = 1$.  Since $g$ acts on $V$ non-trivially, the action on $B$ is non-trivial and $A \otimes B$ is the direct sum of $a$ copies of $B$.

\hspace{1in}

Now take $A = M_1 \pi$ and $B = Y \pi$ so that $V = A \oplus B$.  Suppose that $g$ acts on $A$ as a scalar $c$.  Since $v \ge 7$ and $\dim{B} \le \dim{Y} \le 3$, we see that $a \ge 4$, and by the claim, $c = 1$ and $A \otimes B$ is the direct sum of $a$ copies of $B$.  If $B$ is the zero module, this contradicts the fact that $g$ acts non-trivially on $V$.  Otherwise, $v-a > 0$.  Since $B$ is the image of $Y$, it follows that $A \otimes B \subseteq Y$, and $a(v-a) \le \dim{Y} \le 3$, which is false.  Therefore $g$ does not act on $M_1 \pi$ as a scalar, and hence does not act on $M_1$ as a scalar.  

We may assume that $M_1 = M_f$ where $f$ has degree $u$ over $\F_p$ and $M_1$ and $M_1 \pi$ have types $\alpha'$ and $\beta'$ respectively, where $\beta \subseteq \alpha$.  Write $\alpha = (\alpha_1, \dots, \alpha_s)$ and $\beta = (\beta_1, \dots, \beta_r)$.

\hspace{1in}

\noindent
\emph{Case 1: $u > 1$}.

Writing $S_{M_1}$ for the number of submodules of $M_1$, we have
\begin{eqnarray*}
S_{M_1} &\le& \mathcal{G}_{m_1/u}(q) \\
&\le& C(q) D(q) q^{m_1^2/4u^2} \\
&\le& C(p) D(p) p^{m_1^2/4u} \\
&\le& C(p) D(p) p^{m_1^2/8}.
\end{eqnarray*}
Then 
\begin{eqnarray*}
S_{M} &\le& S_{M_1} \mathcal{G}_y(p) \\
&\le& C(p)^2 D(p)^2 p^{m_1^2/8 + y^2/4} \\
&\le& C(p)^2 D(p)^2 p^{m^2/8 + 9/4} \\
&\le& C(p)^2 D(p)^2 p^{(m-4)^2/4 + 9/4},
\end{eqnarray*}
where the last line uses the fact that $m \ge 14$.  Thus $\log_p{S_M} \le C + (m-4)^2/4$.  

\hspace{1in}

\noindent
\emph{Case 2: $u = 1$}.  

In this case, $f = t-c$ for some $c \in \F_p$.  Since $g$ does not act as a scalar on $M_1$ or $M_1 \pi$, $\alpha_2 \ge \beta_2 > 0$.

By Equation~\ref{logbound},
\[
\log_p{S_M} \le (\alpha_1^2 + \cdots + \alpha_s^2)/4 + s\ep,
\]
so
\[
\log_p{S_M} \le \log_p{S_{M_1}} + \log_p{\mathcal{G}_y(p)} \le (\alpha_1^2 + \cdots + \alpha_s^2 + y^2)/4 + (s+1)\ep.
\]

\hspace{1in}

\noindent
\emph{Subcase 2.1: $\alpha_1 \le m-4$}

If $s = 2$, then 
\begin{eqnarray*}
\log_p{S_M} &\le& (\alpha_1^2 + \alpha_2^2 + y^2)/4 + 3\ep\\
&\le& ((m-4)^2 + 4^2 + 0^2)/4 + 3\ep\\
&\le& (m-4)^2/4 + C.
\end{eqnarray*}
If $s = 3$, then 
\begin{eqnarray*}
\log_p{S_M} &\le& (\alpha_1^2 + \alpha_2^2 + \alpha_3^2 + y^2)/4 + 4\ep \\
&\le& ((m-4)^2 + 3^2 + 1^2 + 0^2)/4 + 4\ep \\
&\le& (m-4)^2/4 + C.
\end{eqnarray*}
Finally, if $4 \le s \le m$, then by Lemma~\ref{quadbound}, we get
\begin{eqnarray*}
\log_p{S_M} &\le& ((m-s)^2 + s)/4 + (s+1)\ep.
\end{eqnarray*}
The right-hand side is maximized at $s = 4$ or $s = m$.  Since $m > 45$ and $\ep \le 6$, it turns out that it is maximized at $s = 4$, where we get a bound of $(m-4)^2/4 + 5\ep + 1$.

\hspace{1in}

\noindent
\emph{Subcase 2.2: $\alpha_1 \ge m-3$}.

So we may assume that $\alpha_1 \ge m-3$.  Then $\alpha_2+\cdots+\alpha_s+y \le 3$, and so $\beta_2 + \cdots + \beta_r + \dim(\pi Y) \le 3$.  Since $\beta_1 + \cdots + \beta_r + \dim(\pi Y) = v \ge 10$, it follows that $\beta_1 \ge 7$ and $\beta_1 - \beta_2 \ge 4$.  Note that $\beta_1 - \beta_2$ is the number of summands of $M_1 \pi$ that are isomorphic to $\F_p[t]/(f-c)$.  So write $M_1 \pi = A \oplus C$, where $a = \dim{A} = \beta_1 - \beta_2$ and $g$ acts as the scalar $c$ on $A$ and not on $C$.  Set $B = C \oplus Y \pi$.  Then $V = A \oplus B$ and by the claim, $c = 1$ and $A \otimes B$ is a direct sum of $a$ copies of $B$.  Then $A \otimes B$ is contained in $Y$ plus the components of $M_1$ that $g$ does not act as a scalar on, so that $a \beta_2 \le \dim{(A \otimes B)} \le \alpha_2 + y \le 3$, which is impossible.
\end{proof}

We can now prove Theorem~\ref{est2thm}, restated here for convenience.

\renewcommand{\arraystretch}{1.2}
\begin{thmx}[Theorem~\ref{est2thm}]
Fix a prime $p$ and integers $d$ and $n$ so that either $n = 2$ and $d \ge 10$ or $n \ge 3$ and $d \ge 3$.  Let $F$ be the free group on $d$ generators and let $d_n$ be the rank of $F_n/F_{n+1}$.  Let
\[
K = \left\{
	\begin{array}{r@{\quad:\quad}l}
		C(p)^5 D(p)^4 p^{17/4} & \textrm{$n = 2$ and $d \ge 10$} \\
		C(p)^2 D(p) p^{3/4} & n \ge 3.
	\end{array}
	\right.
\]
Let
\[
x = \left\{
	\begin{array}{r@{\quad:\quad}l}
		-d & n = 2 \\
		d^2-d_n/2 & n \ge 3.
	\end{array}
	\right.
\]
Then
\begin{enumerate}
\renewcommand{\labelenumi}{\emph{(\alph{enumi})}}
\renewcommand{\theenumi}{\ref{est2thm}(\alph{enumi})}
\item
\[
1 \le \frac{|\mathfrak{C}_{d, n}| \cdot |\GL(d,\F_p)|}{|\mathcal{C}_{d, n}|} \le 1 + K p^x.
\]
\item
\[
1 \le \frac{|\mathfrak{C}_{d, n}|}{|\mathfrak{D}_{d, n}|} \le \frac{1+K p^x}{1-K p^x}.
\]
\end{enumerate}
\end{thmx}
\renewcommand{\arraystretch}{1}

\begin{proof}
Recall that $\mathfrak{C}_{d,n}$ is the set of $\GL(d, \F_p)$-orbits in $\mathcal{C}_{d,n}$, $\mathfrak{D}_{d,n}$ is the set of regular orbits in $\mathfrak{C}_{d,n}$ (that is, the orbits in which every point has trivial stabilizer), and $|\mathcal{C}_{d,n}| = \mathcal{G}_{d_n}(p)$.  If $g \in \GL(d, \F_p)$, then $|(\mathcal{C}_{d,n})^g|$, the number of elements of $\mathcal{C}_{d,n}$ fixed by $g$, is just the number of submodules of $F_n/F_{n+1}$ viewed as a $\F_p \left< g \right>$-module, which we estimated in Theorems~\ref{upperbound} and~\ref{strongerbound}.

We explain first why only the identity element of $\GL(d, \F_p)$ can act as a scalar on $F_n/F_{n+1}$.  By Theorem~\ref{structurethm}, $F_n/F_{n+1}$ has a $\F_p \GL(d, \F_p)$-submodule $M$ which is isomorphic to an extension of $V \wedge V$ by $V$, where $V$ is the natural $\F_p \GL(d,\F_p)$-module.  If $g \in \GL(d,\F_p)$ acts on $F_n/F_{n+1}$ as a scalar $c \in \F_p$, then it acts on $V$ as the scalar $c$, and hence on $V \wedge V$ as the scalar $c^2$.  Thus $c = c^2$ and $c = 1$, so that $g$ is the identity on $V$, that is, the identity element in $\GL(d, \F_p)$.

Suppose first that $n > 2$.  We know from Theorem~\ref{upperbound} that if $g \neq 1$,
\[
|(\mathcal{C}_{d,n})^g| \le C(p)^2 D(p)^2 p^{(d_n^2 - 2d_n + 2)/4}.
\]
By the Cauchy-Frobenius Lemma, 
\begin{eqnarray*}
|\GL(d,\F_p)| \cdot |\mathfrak{C}_{d,n}| &=& \sum_{g \in \GL(d,\F_p)}{|(\mathcal{C}_{d,n})^g|} \\
&=& |\mathcal{C}_{d,n}| + \sum_{g \neq 1}{|(\mathcal{C}_{d,n})^g|} \\
&\le& |\mathcal{C}_{d,n}| + (|\GL(d,\F_p)| - 1) C(p)^2 D(p)^2 p^{(d_n^2-2d_n+2)/4}.
\end{eqnarray*}
By Lemma~\ref{qests} Equation~\ref{gnbounds} and the fact that $2-9p^{(1-d_n)/2}/2 > 1$,
\[
|\mathcal{C}_{d,n}| \ge D(p) p^{d_n^2/4-1/4}.
\]
Since $|\GL(d,\F_p)| \le p^{d^2}$, it follows that
\begin{eqnarray*}
1 &\le& \frac{|\GL(d,\F_p)| \cdot |\mathfrak{C}_{d,n}|}{|\mathcal{C}_{d,n}|} \\
&\le& 1 + C(p)^2 D(p) \; p^{(d_n^2-2d_n+2)/4 + d^2 - d_n^2/4 + 1/4} \\
&=& 1 + K p^{d^2 - d_n/2}.
\end{eqnarray*} 
If $n = 2$, then $F_2/F_3$ is an extension of $V \wedge V$ by $V$, and using the estimates of Lemma~\ref{strongerbound} and the argument above we obtain 
\begin{eqnarray*}
1 &\le& \frac{|\GL(d,\F_p)| \cdot |\mathfrak{C}_{d,n}|}{|\mathcal{C}_{d,n}|} \\
&\le& 1 + K p^{-d}.
\end{eqnarray*}
This proves part $(a)$.

To prove part $(b)$, we observe that $|\mathcal{C}_{d,n}| = \sum{|\GL(d,\F_p)|/|\GL(d,\F_p)_{(w)}|}$, where the sum is over all $\GL(d,\F_p)$-orbits in $\mathcal{C}_{d,n}$ and $|\GL(d,\F_p)_{(w)}|$ is the order of the stabilizer in $\GL(d,\F_p)$ of a typical element $w$ of the orbit under consideration.  Now $|\mathfrak{D}_{d, n}|$ is just the number of orbits for which $|\GL(d,\F_p)_{(w)}| = 1$, so
\[
|\mathcal{C}_{d,n}| \le |\GL(d,\F_p)| \cdot |\mathfrak{D}_{d,n}| + |\GL(d,\F_p)| (|\mathfrak{C}_{d,n}| - |\mathfrak{D}_{d,n}|) / 2.
\]
That is,
\[
(2/|\GL(d,\F_p)|) |\mathcal{C}_{d,n}| - |\mathfrak{C}_{d,n}| \le |\mathfrak{D}_{d,n}|,
\]
so that
\begin{eqnarray*}
\frac{|\mathfrak{C}_{d,n}|}{|\mathfrak{D}_{d,n}|} &\le& \frac{|\mathfrak{C}_{d,n}|}{2|\mathcal{C}_{d,n}|/|\GL(d,\F_p)| - |\mathfrak{C}_{d,n}|} \\
&\le& \frac{|\mathfrak{C}_{d,n}| \cdot |\GL(d,\F_p)| /|\mathcal{C}_{d,n}|}{2 - |\mathfrak{C}_{d,n}| \cdot |\GL(d,\F_p)| /|\mathcal{C}_{d,n}|} \\
&\le& \frac{1 + K p^x}{1 - K p^x}.
\end{eqnarray*}
\end{proof}

\section{Summary} \label{summarysec}

In this section we use Theorems~\ref{bijthm}, \ref{est1thm}, and \ref{est2thm} to prove Theorem~\ref{mainthm} along with two corollaries.

\begin{thmx}[Theorem~\ref{mainthm}]
Fix a prime $p$ and positive integers $d$ and $n$.  Let $r_{d, n}$ be the proportion of $p$-groups minimally generated by $d$ elements and with lower $p$-length at most $n$ whose automorphism group is a $p$-group.  If $n \ge 2$, then
\[
\lim_{d \to \infty}{r_{d, n}} = 1.
\]
If $d \ge 5$, then
\[
\lim_{n \to \infty}{r_{d,n}} = 1.
\]
If 
\begin{eqnarray}
\textrm{$n = 2$ and $d \ge 10$, or $n \ge 3$ and $d \ge 6$, or $n \ge 10$ and $d \ge 5$,} \label{dnconds}
\end{eqnarray}
then
\[
\lim_{p \to \infty}{r_{d,n}} = 1.
\]
\end{thmx}

\begin{proof}
The set of $p$-groups minimally generated by $d$ elements and with lower $p$-length at most $n$ is $\mathfrak{A}_{d,n}$.  When $n = 2$, $\mathfrak{A}_{d,n} = \mathfrak{C}_{d,n}$.  The expression
\[
C(p)^{n-1} D(p)^{n-2} p^{d_{n-1} - d_n/4 + 1/4 + d^2}
\]
goes to $0$ as $d \to \infty$ if $n \ge 3$ or as $n \to \infty$ if $d \ge5$.  If $d$ and $n$ satisfy one of the conditions of Equation~\ref{dnconds}, then the exponent of $p$ is negative.  By Theorem~\ref{est1thm}, it follows that
\begin{eqnarray*}
\lim_{d \to \infty}{\frac{|\mathfrak{A}_{d,n}|}{|\mathfrak{C}_{d,n}|}} &=& 1 \quad \textrm{if $n \ge 2$,} \\
\lim_{n \to \infty}{\frac{|\mathfrak{A}_{d,n}|}{|\mathfrak{C}_{d,n}|}} &=& 1 \quad \textrm{if $d \ge 5$, and} \\
\lim_{p \to \infty}{\frac{|\mathfrak{A}_{d,n}|}{|\mathfrak{C}_{d,n}|}} &=& 1 
\quad \textrm{if one of the conditions in Equation~\ref{dnconds} holds.}
\end{eqnarray*}
The set $\mathfrak{D}_{d,n} \subseteq \mathfrak{C}_{d,n}$ is contained in the subset of $\mathfrak{A}_{d,n}$ of $p$-groups whose automorphism group is a $p$-group.  By Theorem~\ref{est2thm}(b), 
\begin{eqnarray*}
\lim_{d \to \infty}{\frac{|\mathfrak{C}_{d,n}|}{|\mathfrak{D}_{d,n}|}} &=& 1 \quad \textrm{if $n \ge 2$,} \\
\lim_{n \to \infty}{\frac{|\mathfrak{C}_{d,n}|}{|\mathfrak{D}_{d,n}|}} &=& 1 \quad \textrm{if $d \ge 5$, and} \\
\lim_{p \to \infty}{\frac{|\mathfrak{C}_{d,n}|}{|\mathfrak{D}_{d,n}|}} &=& 1 
\quad \textrm{if one of the conditions in Equation~\ref{dnconds} holds.}
\end{eqnarray*}
It follows that $|\mathfrak{A}_{d,n}|/|\mathfrak{D}_{d,n}|$ goes to 1 under the specified limits, and the theorem follows.
\end{proof}

\begin{cor}
Fix a prime $p$ and $n \ge 2$.  Let $s_{d, n}$ be the proportion of $p$-groups generated by at most $d$ elements and with lower $p$-length at most $n$ whose automorphism group is a $p$-group.  Then
\[
\lim_{d \to \infty}{s_{d, n}} = 1.
\]
\end{cor}

\begin{proof}
This follows directly from Theorem~\ref{mainthm} and the trivial observation that the number of $p$-groups generated by at most $d$ elements and with lower $p$-length at most $n$ is finite, while the number of $p$-groups with lower $p$-length at most $n$ is infinite.
\end{proof}

\begin{cor} \label{maincor}
Fix a prime $p$ and $n \ge 2$.  Let $t_{d,n}$ be the proportion of $p$-groups minimally generated by $d$ elements and with lower $p$-length $n$ whose automorphism group is a $p$-group.  Then
\[
\lim_{d \to \infty}{t_{d,n}} = 1.
\]
\end{cor}

\begin{proof}
As $\mathfrak{D}_{d,n} \subseteq \mathfrak{B}_{d,n} \cup \{F_n/F_{n+1}\} \subseteq \mathfrak{A}_{d,n}$, it follows from Theorem~\ref{mainthm} that
\[
\lim_{d \to \infty}{\frac{|\mathfrak{B}_{d,n}|+1}{|\mathfrak{D}_{d,n}|}} = 1.
\]
Since $|\mathfrak{A}_{d,n}| \to \infty$ as $d \to \infty$, Theorem~\ref{mainthm} implies that $|\mathfrak{D}_{d,n}| \to \infty$ as $d \to \infty$, proving that
\[
\lim_{d \to \infty}{\frac{|\mathfrak{B}_{d,n}|}{|\mathfrak{D}_{d,n}|}} = 1.
\]
\end{proof}

Using Theorem~\ref{mainthm}, Henn and Priddy~\cite{hp} prove the following theorem.

\begin{thm}[Henn and Priddy \cite{hp}] \label{hpthm}
Fix a prime $p$ and integers $d, n \ge 2$.  Let $u_{d,n}$ be the proportion of $p$-groups $P$ generated by at most $d$ elements and with lower $p$-length at most $n$ that satisfy the following property: if $H$ is a finite group with Sylow $p$-subgroup $P$, then $H$ has a normal $p$-complement.  Then $\lim_{d \to \infty}{u_{d,n}} = 1$.
\end{thm}

As mentioned in the introduction, the following question remains unanswered.

\begin{question}
Fix a prime $p$.  Let $v_n$ be the proportion of $p$-groups with order at most
$p^n$ whose automorphism group is a $p$-group.  Is it true that $\lim_{n \to \infty}{v_n} = 1$?
\end{question}

\section{Acknowledgements}

We would like to thank Persi Diaconis for introducing us to each other and for his continued support of this project.  We would also like to thank Charles Leedham-Green for several illuminating conversations and for his help with the examples in the introduction.  Finally, we would like to thank Eamonn O'Brien for his help with references and computational data.  For part of this research, the first author was supported by a Department of Defense National Defense Science and Engineering Graduate Fellowship.

\bibliographystyle{amsplain}
\bibliography{grouptheory}

\end{document}